\documentclass[11pt]{amsart}
\usepackage{fullpage}
\usepackage{graphicx}
\usepackage{subfigure}
\usepackage{psfrag}
\usepackage{color}
\usepackage{amsmath}
\usepackage{amssymb}

\newtheorem{theorem}{Theorem}

\newtheorem{remark}[theorem]{Remark}

\title{On the stability of bubble functions and a stabilized mixed finite element formulation for the Stokes problem}

\author{D. Z. Turner} \author{K.~B.~Nakshatrala} \author{K.~D.~Hjelmstad}
\address{Correspondence to: Daniel Z. Turner, Department of Civil and Environmental Engineering, 
  2103 Newmark Laboratory, University of Illinois at Urbana-Champaign, Urbana, Illinois - 61801.}
\email{dzturne1@illinois.edu}
\address{Dr. Kalyana Babu Naskshatrala, Department of Civil and Environmental Engineering, 
  2524 Hydrosystems Laboratory, University of Illinois at Urbana-Champaign, Urbana, Illinois - 61801.}
\email{nakshatr@illinois.edu}
\address{Professor Keith D Hjelmstad, Department of Civil and Environmental Engineering, 
  3129e Newmark Laboratory, University of Illinois at Urbana-Champaign, Urbana, Illinois - 61801.}
\email{kdh@illinois.edu}

\date{\today}

\begin{document}

\begin{abstract}
In this paper we investigate the relationship between stabilized and enriched finite element formulations for the Stokes problem.  We also present a new stabilized mixed formulation for which the stability parameter is derived purely by the method of weighted residuals.  This new formulation allows equal order interpolation for the velocity and pressure fields.  Finally, we show by  counterexample that a direct equivalence between subgrid-based stabilized finite element methods and Galerkin methods enriched by bubble functions cannot be constructed for quadrilateral and hexahedral elements using standard bubble functions.
\end{abstract}


\maketitle



\section{INTRODUCTION}

It is well known that under the mixed Galerkin formulation for the Stokes problem many practically convenient combinations of interpolation functions for the velocity 
and pressure fields often do not yield stable results. In particular, equal order interpolation for velocity and pressure (which is computationally the most convenient) is not stable. This numerical instability is attributed to the lack of stability in the pressure field which is mathematically explained by the celebrated Ladyzhenskaya-Babu\v ska-Brezzi (LBB) stability 
condition \cite{Brezzi}.  One can verify the LBB condition numerically by means of a well designed patch test.

To address the deficiencies in the classical mixed formulation of the Stokes equations, two classes have emerged grouped by similar methodologies: \emph{stabilized finite element methods} and \emph{enriched finite element methods}.  By \emph{stabilized finite element methods} we mean methods that add mesh dependent terms to the standard Galerkin formulation that enable the formulation to satisfy or circumvent the LBB condition \cite{Franca}.  In contrast, \emph{enriched finite element methods} add bubble functions to the finite element function space, which in turn play a stabilizing role.  

Traditionally, the two most popular stabilized methods are the \emph{Streamline-Upwind/Petrov-Galerkin} (SUPG) method \cite{SUPG} and the \emph{Galerkin/least-squares} (GLS) method \cite{GLS}. In the GLS method, least-squares forms of the residuals are added to the Galerkin finite element formulation. These residual based terms are defined over the element interiors only, and the terms on the element boundaries are excluded. The underlying philosophy of the SUPG and the GLS methods is to strengthen the classical variational formulations so that the discrete approximations, which would otherwise be unstable, become stable and convergent. In the mid-90s Hughes revisited the origins of the stabilization schemes from a variational multiscale view point \cite{Hughes2}. Under this variational multiscale method, different stabilization techniques (including GLS and SUPG) are special cases of the underlying subgrid scale modeling concept.  Although referred to as the Douglas-Wang method \cite{Douglas} rather than the variational multiscale method, Franca \textit{et al}. \cite{Franca} proved that the variational multiscale approach for Stokes flow is stable for many elements including T3, TET4, Q4, B8, and combinations of T3 and Q4 elements.

Whereas the applicability of stabilized methods has been established for a wide range of elements, triangular elements have been the primary focus of enriched methods.  Arnold \textit{et al}. \cite{Arnold} were the first to develop the enriched finite element method for Stokes flow.  Using continuous piecewise linear functions enriched by bubbles for velocity and piecewise linear functions for pressure, they showed that the MINI element satisfies the LBB condition.  The stability of similar enrichment for quadrilateral elements has not yet been established.  

Although much work has been reported in the literature (see e.g., \cite{Baiocchi, Brezzi2, Brezzi3, Russo}) showing an equivalence between stabilized and enriched methods, the equivalence is not true for all PDE's nor is it true for all elements.  For example, in the case of nearly incompressible elasticity the equivalence holds for T3, TET4, Q4, and B8 elements (see \cite{Nakshatrala}).  Likewise, for Stokes flow the equivalence holds for T3 and TET4 elements.  However, we will demonstrate that the equivalence breaks down for Stokes flow for Q4 and B8 elements.


The primary aim of this work is twofold: first, we present a new stabilized formulation for the Stokes problem, which can be derived purely by the method of weighted residuals.  This new formulation is stable for equal order interpolation for the velocity and pressure fields, which is computationally convenient.  Second, we show that enriching Q4 and B8 elements with standard bubble functions produces spurious pressure oscillations. The results confirm that one cannot construct an equivalence between stabilized methods and bubble enrichment methods for these elements.

The rest of this paper is organized as follows.  First we present a consistent stabilized formulation for which the stability parameter is constructed from the element residual. Next, we propose an alternative residual based formulation that can be derived purely by the method of weighted residuals.  Lastly, we present a mathematically equivalent enriched formulation and show that whereas the stabilized formulation does not show spurious pressure oscillations for a given test problem, the enriched formulation does.  We conclude with some remarks regarding the equivalence between stabilized and enriched finite element methods for the Stokes problem.

%
  \section{GOVERNING EQUATIONS FOR THE STOKES PROBLEM}
  Let $\Omega$ be a bounded open domain, and $\Gamma$ be its boundary 
  which is assumed to be piecewise smooth. Mathematically, $\Gamma$ 
  is defined as $\Gamma := \bar{\Omega} - \Omega$, where $\bar{\Omega}$ 
  is the closure of $\Omega$. Let the velocity vector field be denoted by 
  $\boldsymbol{v} : \Omega \rightarrow \mathbb{R}^{nd}$, where ``$nd$'' 
  is the number of spatial dimensions. Let the (kinematic) pressure field 
  be denoted by $p:\Omega \rightarrow \mathbb{R}$. 
  As usual, $\Gamma$ is divided into two parts, denoted by 
  $\Gamma^{\boldsymbol{v}}$ and $\Gamma^{\boldsymbol{t}}$, such 
  that $\Gamma^{\boldsymbol{v}} \cap \Gamma^{\boldsymbol{t}} =\emptyset$ 
  and $\Gamma^{\boldsymbol{v}} \cup \Gamma^{\boldsymbol{t}}=\Gamma$.  $\Gamma^{\boldsymbol{v}}$ is the part of the boundary 
  on which velocity is prescribed, and $\Gamma^{\boldsymbol{t}}$ 
  is part of the boundary on which traction is prescribed. 
  The governing equations for Stokes flow can be written as 
  \begin{align}
    \label{Eqn:SNS_Equilibrium}
    -2\nu \nabla^2 \boldsymbol{v} + 
    \nabla p &= \boldsymbol{b} \qquad \; \;  \ \mbox{in}  \quad \Omega \\
    \label{Eqn:SNS_Continuity}
    \nabla \cdot \boldsymbol{v} &= 0 \qquad \; \;   \ \mbox{in} \quad \Omega \\
    \label{Eqn:SNS_VelocityBC}
    \boldsymbol{v} &= \boldsymbol{v}^{\mathrm{p}} 
    \qquad \; \mbox{on} \quad \Gamma^{\boldsymbol{v}} \\
    \label{Eqn:SNS_TractionBC}
    -p \boldsymbol{n} + \nu (\boldsymbol{n} \cdot \nabla) \boldsymbol{v} &= 
    \boldsymbol{t}^{\boldsymbol{n}} \qquad \; \mbox{on} \quad \Gamma^{\boldsymbol{t}} 
  \end{align}
  where $\nabla$ is the gradient operator, $\nabla^2$ is the laplacian operator, $\boldsymbol{b}$ is the body force, 
  $\nu > 0$ is the kinematic viscosity, $\boldsymbol{v}^{\mathrm{p}}$ is the 
  prescribed velocity vector field, $\boldsymbol{t}^{\boldsymbol{n}}$ is the 
  prescribed traction, and $\boldsymbol{n}$ is the unit outward normal vector 
  to $\Gamma$. Equation \eqref{Eqn:SNS_Equilibrium} represents the balance of 
  linear momentum, and equation \eqref{Eqn:SNS_Continuity} represents the continuity 
  equation for an incompressible continuum. Equations \eqref{Eqn:SNS_VelocityBC} 
  and \eqref{Eqn:SNS_TractionBC} are the Dirichlet and Neumann boundary 
  conditions, respectively. 
  

  In the next section, we present the classical mixed formulation for 
  the Stokes equations which will be the basis for the stabilized and enriched formulations. 
  
%
\section{CLASSICAL MIXED FORMULATION}
Before we present the classical mixed formulation for the Stokes equations, 
let us define the function spaces that will be used in the remainder of this paper. 
The function spaces for the velocity $\boldsymbol{v}(\boldsymbol{x})$ and the weighting 
function associated with velocity, denoted by $\boldsymbol{w}(\boldsymbol{x})$, are 
respectively defined as
\begin{align}
  \label{Eqn:SNS_Function_Space_for_v}
  \mathcal{V} &:= \{\boldsymbol{v}  \; \big| \; \boldsymbol{v} \in (H^{1}(\Omega))^{nd},
  \boldsymbol{v} = \boldsymbol{v}^{\mathrm{p}} \; \mathrm{on} \; \Gamma^{\boldsymbol{v}} \} \\
  \label{Eqn:SNS_Function_Space_for_w}
  \mathcal{W} &:= \{\boldsymbol{w}  \; \big| \; \boldsymbol{w} \in (H^{1}(\Omega))^{nd},
  \boldsymbol{w} = \boldsymbol{0} \; \mathrm{on} \; \Gamma^{\boldsymbol{v}} \}
\end{align}
where $H^1(\Omega)$ is a standard Sobolev space \cite{Brezzi}. In the classical mixed 
formulation the function space for the pressure $p(\boldsymbol{x})$ and its corresponding 
weighting function $q(\boldsymbol{x})$ are given by 
%
\begin{equation}
  \label{Eqn:SNS_Function_Space_for_p}
  \mathcal{P} := \{p \; \big| \; p \in L^{2}(\Omega) \}
\end{equation}
where $L^2(\Omega)$ is the space of square-integrable functions on the domain $\Omega$. In the stabilized formulations, the function space for $p(\boldsymbol{x})$ and $q(\boldsymbol{x})$ will be defined as
 \begin{equation}
  \label{Eqn:SNS_HVM_Function_Space_for_p}
  \underline{\mathcal{P}} := \{p \; \big| \; p \in H^{1}(\Omega) \}
\end{equation}
For further details on function spaces refer to Brezzi and Fortin \cite{Brezzi}. 
%
\begin{remark}
  When Dirichlet boundary conditions are imposed everywhere on the boundary, 
  that is $\Gamma^{\boldsymbol{t}} = \emptyset$, the pressure 
  can be determined only up to an arbitrary constant. In order to define the 
  pressure field uniquely, it is common to prescribe the average value of pressure,
  \begin{align}
    \int_{\Omega} p \; \mathrm{d} \Omega = p_0 
  \end{align}
  where $p_0$ is arbitrarily chosen (and can be zero).  Then, the appropriate function spaces for the pressure that should be used instead of 
  $\mathcal{P}$ (defined in equation \eqref{Eqn:SNS_Function_Space_for_p}) is 
  \begin{align}
    \mathcal{P}_0 := \{p \; \big| \; p \in L^{2}(\Omega), \int_{\Omega} p \; \mathrm{d} \Omega  = 0\}
  \end{align}
  Another way to define the pressure uniquely is to prescribe the value of the 
  pressure at a point, which is computationally the most convenient. 
\end{remark}

The classical mixed formulation (which is based on the Galerkin principle) for the 
Stokes equations can be written as: Find $\boldsymbol{v}(\boldsymbol{x}) 
\in \mathcal{V}$ and $p(\boldsymbol{x}) \in \mathcal{P}$ such that 
\begin{align}
  \label{Eqn:SNS_CM_Momentum}
  a(\boldsymbol{w};\boldsymbol{v}) + 
  b(\boldsymbol{w};p) = 
  f(\boldsymbol{w}) & \quad \forall \ \boldsymbol{w} \in \mathcal{W} \\
  \label{Eqn:SNS_CM_Continuity}
  b(\boldsymbol{v};q) = 0 & \quad \forall \ q \in \mathcal{P}
\end{align}
Let us define the bilinear forms as: 
\begin{align}
  \label{Eqn:SNS_a_Galerkin}
  a(\boldsymbol{w};\boldsymbol{v}) &:= 
  \int_{\Omega} \nabla \boldsymbol{w}:2\nu \nabla 
  \boldsymbol{v} \; \mathrm{d} \Omega \\
  \label{Eqn:SNS_b_Galerkin}
  b(\boldsymbol{w};p) &:= 
  -\int_{\Omega} \left(\nabla \cdot \boldsymbol{w}\right)
  \; p \; \mathrm{d} \Omega 
\end{align}
and the linear functional as 
\begin{align}
  \label{Eqn:SNS_f_Galerkin}
  f(\boldsymbol{w}) & := \int_{\Omega} \boldsymbol{w} \cdot 
  \boldsymbol{b} \; \mathrm{d} \Omega + \int_{\Gamma^{\boldsymbol{t}}} 
  \boldsymbol{w} \cdot \boldsymbol{t}^{\boldsymbol{n}} \; \mathrm{d} \Gamma
\end{align}

Once the weak formulation of the governing equations is established, the approximate 
solution based on the finite element method is determined in the usual manner. First 
one chooses the approximating finite element spaces, which (for a conforming formulation) 
will be finite dimensional subspaces of the underlying function spaces of the weak 
formulation. Let the finite element function spaces for the velocity, the weighting 
function associated with the velocity, and the pressure be denoted by $\mathcal{V}^{h} 
\subseteq \mathcal{V}$, $\mathcal{W}^h \subseteq \mathcal{W}$, and $\mathcal{P}^h 
\subseteq \mathcal{P}$ respectively. The finite element formulation of the classical 
mixed formulation reads: Find $\boldsymbol{v}^h(\boldsymbol{x}) \in \mathcal{V}^h$ and 
$p^h(\boldsymbol{x}) \in \mathcal{P}^h$ such that 
\begin{align}
  \label{Eqn:SNS_FEM_CM_Momentum}
  a(\boldsymbol{w}^h;\boldsymbol{v}^h) + 
  b(\boldsymbol{w}^h;p^h) = 
  f(\boldsymbol{w}^h) & \quad \forall \ \boldsymbol{w}^h \in \mathcal{W}^h \\
  \label{Eqn:SNS_FEM_CM_Continuity}
  b(\boldsymbol{v}^h;q^h) = 0 & \quad \forall \ q^h \in \mathcal{P}^h
\end{align}
For mixed formulations, the inclusions $\mathcal{V}^h \subseteq \mathcal{V}$, $\mathcal{W}^h 
\subseteq \mathcal{W}$, and $\mathcal{P}^h \subseteq \mathcal{P}$ are 
themselves not sufficient to produce stable results, and additional 
conditions must be met by these finite element spaces to obtain 
meaningful numerical results. A systematic study of these types of 
conditions on function spaces to obtain stable numerical results 
is the main theme of \emph{mixed finite elements}. One of the main 
conditions to be met is the LBB \emph{inf-sup} stability condition. 

Although the classical mixed formulation has many advantages (mainly its simplicity and extensions to turbulent flows), it also has several numerical 
deficiencies.  Most importantly, many combinations of shape functions for the 
velocity and pressure do not satisfy the LBB 
stability condition and therefore exhibit unphysical oscillations in numerical simulations. 
As mentioned previously, two classes of methods have been developed to overcome the
limitations associated with the classical Galerkin approach; methods that augment the formulation with stabilizing terms to circumvent the LBB stability condition and those that enrich the function space to satisfy the LBB condition.
%
\section{VARIATIONAL MULTISCALE FRAMEWORK}
Hughes \cite{Hughes2} proposed a variational framework based on the 
multiscale decomposition of the underlying fields into a coarse or resolvable scale and a subgrid or unresolvable scale. This framework provides a systematic 
procedure to develop stable finite element formulations. In this section, we 
present a multiscale formulation for the Stokes 
equations. A similar formulation for Darcy flow is presented in \cite{Nakshatrala2}.
%
\subsection{Multiscale decomposition}
Let us divide the domain $\Omega$ into $N$ non-overlapping subdomains 
$\Omega^e$ (which in the finite element context will be elements) such 
that
%
\begin{equation}
  \Omega = \overset{N}{\underset{e = 1}{\bigcup}} \Omega^{e}
\end{equation}
The boundary of the element $\Omega^{e}$ is denoted by $\Gamma^{e}$. We 
decompose the velocity field $\boldsymbol{v}(\boldsymbol{x})$ into 
coarse-scale and fine-scale components, indicated as $\bar{\boldsymbol{v}}
(\boldsymbol{x})$ and $\boldsymbol{v}'(\boldsymbol{x})$, respectively. 
To wit, 
%
\begin{equation}
  \label{Eqn:SNS_Decompose_v}
  \boldsymbol{v}(\boldsymbol{x}) = \bar{\boldsymbol{v}} 
  (\boldsymbol{x}) + \boldsymbol{v}'(\boldsymbol{x})
\end{equation}
Likewise, we decompose the weighting function $\boldsymbol{w}(\boldsymbol{x})$ 
into coarse-scale $\bar{\boldsymbol{w}}(\boldsymbol{x})$ and fine-scale 
$\boldsymbol{w}'(\boldsymbol{x})$ components.
%
\begin{equation}
  \label{Eqn:SNS_Decompose_w}
  \boldsymbol{w}(\boldsymbol{x}) = \bar{\boldsymbol{w}}(\boldsymbol{x}) 
  + \boldsymbol{w}'(\boldsymbol{x})
\end{equation}
We further make an assumption that the fine-scale 
components vanish along each element boundary.
%
\begin{equation}
  \label{Eqn:SNS_Fine_Scale_Vanish}
  \boldsymbol{v}'(\boldsymbol{x}) = \boldsymbol{w}'(\boldsymbol{x}) 
  = \boldsymbol{0} \quad \mbox{on} \quad \Gamma^{e} \; ; \; e = 1, \ldots, N
\end{equation}
Let $\bar{\mathcal{V}}$ be the function space for the coarse-scale component 
of the velocity $\bar{\boldsymbol{v}}$, and $\bar{\mathcal{W}}$ be the 
function space for $\bar{\boldsymbol{w}}$; and are defined as 
\begin{align}
  \bar{\mathcal{V}} := \mathcal{V}; \;
  \bar{\mathcal{W}} := \mathcal{W}
\end{align}
where $\mathcal{V}$ and $\mathcal{W}$ are defined earlier in equation 
\eqref{Eqn:SNS_Function_Space_for_v} and equation \eqref{Eqn:SNS_Function_Space_for_w} 
respectively. Let $\mathcal{V}'$ be the function space for both the fine-scale 
component of the velocity $\boldsymbol{v}'$ and its corresponding weighting 
function $\boldsymbol{w}'$, and is defined as 
\begin{align}
  \label{Eqn:SNS_function_space_for_v_prime}
  \mathcal{V}' := \{\boldsymbol{v} \; \bigr| \; \boldsymbol{v} 
  \in (H^{1}(\Omega^e))^{nd} , \; \boldsymbol{v} = \boldsymbol{0} \; \mbox{on} \; 
  \Gamma^e, e = 1, \ldots, N \} 
\end{align}
The velocity field $\boldsymbol{v}(\boldsymbol{x})$ is now an element of the function 
space generated by the direct sum of $\bar{\mathcal{V}}$ and $\mathcal{V}'$, denoted 
by $\bar{\mathcal{V}} \oplus \mathcal{V}'$. Similarly the direct sum of $\bar{\mathcal{W}}$ 
and $\mathcal{V}{'}$, denoted by $\bar{\mathcal{W}}\oplus\mathcal{V}'$, is the function 
space for the field $\boldsymbol{w}(\boldsymbol{x})$. 

In theory, we could decompose the pressure field into coarse-scale and fine-scale 
components. However, for simplicity we assume that there are no fine-scale terms 
for the pressure $p(\boldsymbol{x})$ and for its corresponding weighting function 
$q(\boldsymbol{x})$. Hence the function space for the fields $p(\boldsymbol{x})$ 
and $q(\boldsymbol{x})$ is $\mathcal{P}$.
%
\subsection{Two-level classical mixed formulation}
Substitution of equations \eqref{Eqn:SNS_Decompose_v} and \eqref{Eqn:SNS_Decompose_w} 
into the classical mixed formulation given by equations \eqref{Eqn:SNS_CM_Momentum} and 
\eqref{Eqn:SNS_CM_Continuity} becomes the first point of departure from the 
classical Galerkin formulation. 
%
\begin{align}
  \label{Eqn:SNS_Multiscale_Momentum}
  a(\bar{\boldsymbol{w}} + \boldsymbol{w}';\bar{\boldsymbol{v}} + \boldsymbol{v}') + 
  b(\bar{\boldsymbol{w}} + \boldsymbol{w}';p) 
  &= f(\bar{\boldsymbol{w}} + \boldsymbol{w}')  \\
  \label{Eqn:SNS_Multiscale_Continuity}
  b(\bar{\boldsymbol{v}} + \boldsymbol{v}';q) &= 0 
\end{align}

Because the weighting functions $\bar{\boldsymbol{w}}$ and $\boldsymbol{w}'$ 
are arbitrary, and because the functionals are linear in the weighting functions, we can write 
the above problem as two sub-problems. The \emph{coarse-scale problem} can be written as:
\begin{align}
  \label{Eqn:SNS_CoarseScale_Momentum}
  a(\bar{\boldsymbol{w}};\bar{\boldsymbol{v}} + \boldsymbol{v}') + 
  b(\bar{\boldsymbol{w}};p) 
  &= f(\bar{\boldsymbol{w}}) \quad \forall \ \bar{\boldsymbol{w}} \in \bar{\mathcal{W}} \\
  \label{Eqn:SNS_CoarseScale_Continuity}
  b(\bar{\boldsymbol{v}} + \boldsymbol{v}';q) &= 0 \quad \forall \ q \in \mathcal{P}
\end{align}
where the quantities $a(\cdot;\cdot)$, $b(\cdot;\cdot)$ and $f(\cdot;\cdot)$ 
are defined in equations \eqref{Eqn:SNS_a_Galerkin}-\eqref{Eqn:SNS_f_Galerkin}. The 
\emph{fine-scale problem} can be written as:
\begin{align}
  \label{Eqn:SNS_FineScale_Momentum}
  a(\boldsymbol{w}';\bar{\boldsymbol{v}} + \boldsymbol{v}') + 
  b(\boldsymbol{w}';p) 
  &= f(\boldsymbol{w}') \quad \forall \ \boldsymbol{w}' \in \mathcal{W}' 
\end{align}
%
  %
  %
  %
  %
%
\begin{remark}
Note that the fine scale problem is independent and uncoupled at the element level (defined over the sum of element interiors).  Due to the assumption that the subgrid scale response vanishes on the element boundaries, $a(\bar{\boldsymbol{w}};\bar{\boldsymbol{v}} + \boldsymbol{v}') = a(\bar{\boldsymbol{w}};\bar{\boldsymbol{v}}) + a(\bar{\boldsymbol{w}};\boldsymbol{v}')$.
\end{remark}
%
Using the linearity of the solution field and the divergence theorem on $a(\bar{\boldsymbol{w}};\boldsymbol{v}')$, we may alternatively write the \emph{coarse-scale problem} as:
\begin{align}
  \label{Eqn:SNS_CoarseScale_Momentum_Alt}
  a(\bar{\boldsymbol{w}};\bar{\boldsymbol{v}}) + 
  b(\bar{\boldsymbol{w}};p) + 
  c(\bar{\boldsymbol{w}};\boldsymbol{v}') 
  &= f(\bar{\boldsymbol{w}}) \\
  \label{Eqn:SNS_CoarseScale_Continuity_Alt} 
  b(\bar{\boldsymbol{v}};q) +
  d(\boldsymbol{v}';q) &= 0
\end{align}
and the \emph{fine-scale problem} as: 
\begin{align}
  \label{Eqn:SNS_FineScale_Momentum_Alt}
  a(\boldsymbol{w}';\boldsymbol{v}') + 
  c(\bar{\boldsymbol{v}};\boldsymbol{w}') +
  d(\boldsymbol{w}';p)  &= f(\boldsymbol{w}')
\end{align}
where 
\begin{align}
\label{Eqn:SNS_New_Bilinear_Forms}
c(\boldsymbol{w};\boldsymbol{v}) := -\int_{\Omega} \nabla^2 \boldsymbol{w} \cdot 2\nu \boldsymbol{v} \; \mathrm{d} \Omega \\
d(\boldsymbol{w};p) :=  \int_{\Omega} \boldsymbol{w} \cdot \nabla p \; \mathrm{d} \Omega
\end{align}

%
\section{FINE--SCALE INTERPOLATION AND BUBBLE FUNCTIONS}
%
If one chooses a single bubble function for interpolating the fine-scale variables (similar to the MINI element), then we have
\begin{align}
 \label{eq:DefineBubbles}
\boldsymbol{v}' = b^e \boldsymbol{\beta}; \quad \boldsymbol{w}' = b^e \boldsymbol{\gamma}
\end{align}
where $b^e$ is a bubble function, and $\boldsymbol{\beta}$ and $\boldsymbol{\gamma}$ 
are constant vectors. 
The gradients of the fine-scale velocity and weighting functions are
\begin{align}
\nabla \boldsymbol{v}' = \boldsymbol{\beta}\nabla b^{eT}; \quad \nabla \boldsymbol{w}' = \boldsymbol{\gamma}\nabla b^{eT}
\end{align}
where $\nabla b^e $ is a $\mathrm{dim} \times 1 $ vector of the derivatives of the bubble function.
Standard bubble functions for several elements are provided in Table \ref{table:BubbleDefs}. 
\begin{table}[htb]
\caption{Bubble Functions for Standard Finite Elements}
\centering
\begin{tabular}{l l}
\hline
Element   & Bubble function \\ 
\hline
T3 & $\xi_1\xi_2(1 - \xi_1 - \xi_2) $ \\ 
TET4 & $\xi_1\xi_2\xi_3(1 - \xi_1 - \xi_2 - \xi_3) $ \\
Q4 & $ (1-\xi_1^2)(1-\xi_2^2)$\\
B8 & $ (1-\xi_1^2)(1-\xi_2^2)(1-\xi_3^2)$ \\
\hline
\end{tabular}
\label{table:BubbleDefs}
\end{table}

We shall substitute these expressions into the above subproblems in two different fashions, which  brings us to the point of departure between stabilized and enriched methods.
%
\subsection{Weak variational multiscale formulation}
In the spirit of a stabilized method, we eliminate the fine-scale variables by solving the \emph{fine-scale problem} (equation \eqref{Eqn:SNS_FineScale_Momentum_Alt}) in terms of the coarse-scale variables.  We then substitute the fine-scale solution into the 
\emph{coarse-scale problem} (equation \eqref{Eqn:SNS_CoarseScale_Momentum_Alt}) and solve 
the \emph{coarse-scale problem} to obtain $\bar{\boldsymbol{v}}(\boldsymbol{x})$ 
and $p(\boldsymbol{x})$. This procedure also produces the familiar stabilization parameter, $\boldsymbol{\tau}$, with which we augment the classical Galerkin formulation.  Traditionally, one solves the \emph{fine-scale problem} in terms of the coarse-scale variables in a weak or integral sense.  For this reason, we refer to this method as the weak variational multiscale (WVM) formulation. 

\subsection{Stabilization parameter}
Typically, the stabilization parameter is derived in a consistent manner by incorporating the coarse-scale residual evaluated over the element.  Examples of such formulations include the work of Masud and Khurram \cite{Masud} for the Stokes equations and that of Nakshatrala \emph{et al} \cite{Nakshatrala} for nearly incompressible linear elasticity.  The derivation proceeds as follows.

Returning to equation \eqref{Eqn:SNS_FineScale_Momentum_Alt}, substituting equation \eqref{eq:DefineBubbles}, and noting the arbitrariness of $\boldsymbol{\gamma}$ we have
\begin{align}
  \label{Eqn:SNS_Stabilized_Fine_1}
  2\nu \int_{\Omega^e} \left| \nabla b^e \right|^2 \; \mathrm{d} \Omega \boldsymbol{\beta}
  = \int_{\Omega^e} b^e \bar{\boldsymbol{r}} \; \mathrm{d} \Omega
\end{align}
where $\bar{\boldsymbol{r}} := 2\nu \nabla^2 \bar{\boldsymbol{v}} - \nabla p + \boldsymbol{b}$ is the collection of the coarse-scale terms in the \emph{fine-scale problem}. To solve for $\boldsymbol{\beta}$, one can make the approximation that in the limit of mesh refinement, the coarse-scale residual is constant over the element domain.  Hence, $\bar{\boldsymbol{r}}$ is moved outside of the integral in equation \eqref{Eqn:SNS_Stabilized_Fine_1} such that 
\begin{align}
  \label{Eqn:SNS_Stabilized_Fine_2}
  \boldsymbol{\beta} = \left[ 2\nu \int_{\Omega^e} \left| \nabla b^e \right|^2 \; \mathrm{d} \Omega \right]^{-1} \int_{\Omega^e} b^e \; \mathrm{d} \Omega \ \bar{\boldsymbol{r}}
\end{align}
\begin{remark}
Note this is the only approximation introduced for this method, aside from the assumption that the subgrid scales vanish on the boundary (which is the key feature of the variational multiscale framework).
\end{remark}
\begin{remark}
\label{remark:ResidRemark}
In the case of T3 and TET4 elements, the statement that the coarse-scale residual is constant over the element domain is not an approximation, but is exactly true if $\boldsymbol{b}$ is constant.
\end{remark}
%
%
Referring to equation \eqref{eq:DefineBubbles}, the fine-scale velocity may then be written as
\begin{align}
  \label{Eqn:SNS_Define_Velocity_Fine}
  \boldsymbol{v}' = b^e \boldsymbol{\beta} = \frac{1}{2\nu}\boldsymbol{\tau}\bar{\boldsymbol{r}}
\end{align}
where we have introduced the stabilization parameter $\boldsymbol{\tau}$
\begin{align}
  \label{Eqn:SNS_Define_Tau}
  \boldsymbol{\tau} = b^e \left[\int_{\Omega^e}  \left| \nabla b^e \right|^2 \; \mathrm{d} \Omega \right]^{-1} \int_{\Omega^e} b^e \; \mathrm{d} \Omega
\end{align}

\subsection{Weak variational multiscale Galerkin formulation}
Since we have an expression for the fine-scale velocity, we can substitute equation \eqref{Eqn:SNS_Define_Velocity_Fine} back into the the \emph{coarse-scale problem} to obtain a stabilized version of the Galerkin formulation
\begin{align}
  \label{Eqn:SNS_Stabilized_Momentum}
  a(\bar{\boldsymbol{w}};\bar{\boldsymbol{v}}) + 
  b(\bar{\boldsymbol{w}};p) -
  c(\bar{\boldsymbol{w}} ;\frac{1}{2\nu}\boldsymbol{\tau}(2\nu \nabla^2 \bar{\boldsymbol{v}} - \nabla p + \boldsymbol{b})) = 
  f(\bar{\boldsymbol{w}}) & \quad \forall \ \bar{\boldsymbol{w}} \in \bar{\mathcal{W}} \\
  \label{Eqn:SNS_stabilized_Continuity}
  b(\bar{\boldsymbol{v}};q) - d(\frac{1}{2\nu}\boldsymbol{\tau}(2\nu \nabla^2 \bar{\boldsymbol{v}} - \nabla p + \boldsymbol{b});q) = 0 & \quad \forall \ q \in \underline{\mathcal{P}}
\end{align}
Note the bilinear forms are defined in equations \eqref{Eqn:SNS_a_Galerkin}--\eqref{Eqn:SNS_f_Galerkin} and \eqref{Eqn:SNS_New_Bilinear_Forms}.

%
\section{STRONG VARIATIONAL MULTISCALE FORMULATION}
We now present a new stabilized formulation for the Stokes problem that is consistently derived from the method of weighted residuals.  Whereas traditionally the \emph{fine-scale problem} is solved in a weak or integral sense, in the following formulation we solve the \emph{fine-scale problem} in a strong sense.  Therefore, we refer to this method as the strong variational multiscale (SVM) formulation.  

Using integration by parts and the linearity of the solution field, we may rewrite the \emph{fine-scale problem} (given by equation \eqref{Eqn:SNS_FineScale_Momentum_Alt}) as:
\begin{align}
  \label{Eqn:SNS_FineScale_Momentum_Alternative1}
  c(\boldsymbol{v}';\boldsymbol{w}') +
  c(\bar{\boldsymbol{v}};\boldsymbol{w}') +
  b(\boldsymbol{w}';p) &=
  f(\boldsymbol{w}')
\end{align}
Using the notation for the coarse-scale residual $\bar{\boldsymbol{r}} := 2\nu \nabla^2 \bar{\boldsymbol{v}} - \nabla p + \boldsymbol{b}$, the above equation can be written as,
\begin{align}
  \label{Eqn:SNS_FineScale_Momentum_Alternative2}
  \int_{\Omega} \boldsymbol{w}' \cdot (- 2\nu \nabla^2 \boldsymbol{v}' -  \bar{\boldsymbol{r}})\; \mathrm{d} \Omega = 0
 \end{align}
Because $\boldsymbol{w}'$ is arbitrary and vanishes on the element boundaries and because $\boldsymbol{v}'$ is constrained to vanish on the element boundaries, the strong form of equation \eqref{Eqn:SNS_FineScale_Momentum_Alternative2} is
\begin{align}
  \label{Eqn:SNS_FineScale_Momentum_Alt3}
  2\nu \nabla^2 \boldsymbol{v}' &= -  \bar{\boldsymbol{r}} \qquad \; \;  \ \mbox{in}  \quad \Omega_e \; ; e = 1 \ldots N\\
    \boldsymbol{v}' &= 0 \qquad \; \; \; \; \;  \ \mbox{on}  \quad \Gamma_e
 \end{align}
\begin{remark}
The strong form may also be written as 
\begin{align}
  \label{Eqn:SNS_FineScale_Alt_Strong_Form}
  \mathcal{L} \left[ \boldsymbol{v}' \right] = -  \bar{\boldsymbol{r}}(\boldsymbol{x}) \;   \mbox{in} \; \Omega_e \; ; \;  \boldsymbol{v}'(\boldsymbol{x}) = 0 \; \mbox{on} \; \Gamma_e \; ; \; e = 1 \ldots N
\end{align}
where $\mathcal{L} \left[ \cdot \right] = 2\nu \nabla^2(\cdot)$ is the linear differential operator of the \emph{fine-scale problem}.  The analytical solution to equation \eqref{Eqn:SNS_FineScale_Alt_Strong_Form} over the element domain may be written as
 \begin{align}
  \label{Eqn:SNS_Greens_Analytical}
  \boldsymbol{v}'(\boldsymbol{x}) = - \int_{\Omega^e} \boldsymbol{G}(\boldsymbol{x}, \boldsymbol{y})\bar{\boldsymbol{r}}(\boldsymbol{y})  \; \mathrm{d} \Omega_{\boldsymbol{y}}
 \end{align}
 where $\boldsymbol{G}(\boldsymbol{x}, \boldsymbol{y})$ is the Green's function for the operator $\mathcal{L}$.  The potential for $\boldsymbol{\tau}$ to emanate from the element's Green's function has been pointed out in \cite{Hughes2}.
\end{remark}
Obtaining an analytical solution for the Green's function that is valid for any element configuration is not always possible.  Also, in order to get stable results an approximation to the Green's function will suffice.  To this end, we approximate the solution using a single bubble function
\begin{align}
 \label{eq:DefineBubbleVelocity}
\boldsymbol{v}' = b^e \boldsymbol{\beta}; \; \nabla^2 \boldsymbol{v}' = (\nabla^2 b^e) \boldsymbol{\beta}
\end{align}
where $\nabla^2 b^e$ is defined as the Laplacian of the bubble function, which will never be zero (see Appendix).  Substituting equation \eqref{eq:DefineBubbleVelocity} into equation \eqref{Eqn:SNS_FineScale_Momentum_Alt3} we have
\begin{align}
 \label{eq:AlmostTau}
 \boldsymbol{\beta} = - \frac{1}{2\nu \nabla^2 b^e} \bar{\boldsymbol{r}}
\end{align}
We now have an expression for the fine-scale velocity $\boldsymbol{v}'$
\begin{align}
 \label{eq:DefineAlternateFineVel}
 \boldsymbol{v}' = - \frac{b^e}{2\nu \nabla^2 b^e} \bar{\boldsymbol{r}} = - \frac{1}{2\nu} \boldsymbol{\tau}\bar{\boldsymbol{r}}
\end{align}
where $\boldsymbol{\tau}$ is the stabilization parameter defined as
\begin{align}
 \label{eq:DefineTau}
 \boldsymbol{\tau} := \frac{b^e}{\nabla^2 b^e} 
\end{align}
A straightforward analysis shows that for an element with characteristic dimension $h$, the stabilization parameter $\boldsymbol{\tau}$ scales as $h^2$.
\begin{remark}
It is well-known in the mixed finite element literature (for example, see \cite{Douglas,PSPG}) that $\boldsymbol{\tau}$ must scale as $h^2$ to guarantee convergence, which appears to be satisfied by \eqref{eq:DefineTau}.
\end{remark}
%
%
\begin{remark}
The above stabilization parameter makes no approximations relative to the coarse-scale residual, as in the weak variational multiscale formulation.  Therefore in the case of quadrilateral or hexahedral elements, no additional approximations are introduced preserving a mathematically exact correspondence with the enriched formulation presented below.
\end{remark}
\subsection{Weak variational multiscale Galerkin formulation}
After substitution of equation \eqref{eq:DefineAlternateFineVel} into the \emph{coarse-scale problem} (equations \eqref{Eqn:SNS_CoarseScale_Momentum_Alt} and \eqref{Eqn:SNS_CoarseScale_Continuity_Alt})
the resulting weak form is again expressed exactly as equations \eqref{Eqn:SNS_Stabilized_Momentum} and \eqref{Eqn:SNS_stabilized_Continuity}.
\begin{align}
  \label{Eqn:SNS_New_Stabilized_Momentum}
  a(\bar{\boldsymbol{w}};\bar{\boldsymbol{v}}) + 
  b(\bar{\boldsymbol{w}};p) +
  c(\bar{\boldsymbol{w}} ;\frac{1}{2\nu}\boldsymbol{\tau}(2\nu \nabla^2 \bar{\boldsymbol{v}} - \nabla p + \boldsymbol{b}) ) = 
  f(\bar{\boldsymbol{w}}) & \quad \forall \ \bar{\boldsymbol{w}} \in \bar{\mathcal{W}} \\
  \label{Eqn:SNS_New_Stabilized_Continuity}
  b(\bar{\boldsymbol{v}};q) + d(\frac{1}{2\nu}\boldsymbol{\tau}(2\nu \nabla^2 \bar{\boldsymbol{v}} - \nabla p + \boldsymbol{b});q) = 0 & \quad \forall \ q \in \underline{\mathcal{P}}
\end{align}
Note that for linear elements like the T3 and TET4, $\nabla^2 \bar{\boldsymbol{w}}$ and $\nabla^2 \bar{\boldsymbol{v}}$ will be exactly zero.
  %
%
\section{ENRICHED FORMULATION}
For the enriched formulation we treat the \emph{coarse} and \emph{fine-scale problems} (equations \eqref{Eqn:SNS_CoarseScale_Momentum}-- \eqref{Eqn:SNS_FineScale_Momentum}) as two residual equations of the variables $\bar{\boldsymbol{v}}$, $\boldsymbol{v}'$, and $p$. Instead of analytically solving for $\boldsymbol{v}'$ in terms of the coarse-scale variables (as in a stabilized formulation), we use static condensation to solve the problem in a two stage manner.  The emphasis in this section is placed on the solution strategy since it represents the most relevant features of the enriched formulation.  
\subsection{Scalar residual}
The scalar residual equations may be written as
\begin{align}
  \label{Eqn:SNS_Scalar_Resid_Coarse}
  r_c(\bar{\boldsymbol{v}};\boldsymbol{v}',p) &:= a(\bar{\boldsymbol{w}};\bar{\boldsymbol{v}} + \boldsymbol{v}') + 
  b(\bar{\boldsymbol{w}};p) 
  - f(\bar{\boldsymbol{w}}) \\
  \label{Eqn:SNS_Scalar_Resid_Pressure}
  r_p(\bar{\boldsymbol{v}};\boldsymbol{v}') &:= b(\bar{\boldsymbol{v}} + \boldsymbol{v}';q) \\
  \label{Eqn:SNS_Scalar_Resid_Fine}
  r_f(\bar{\boldsymbol{v}};\boldsymbol{v}',p) &:= a(\boldsymbol{w}';\bar{\boldsymbol{v}} + \boldsymbol{v}') + 
  b(\boldsymbol{w}';p) - f(\boldsymbol{w}')
\end{align}
where the subscripts `c', `p', and `f' stand for \emph{coarse}, \emph{pressure}, and \emph{fine}.
\subsection{Vector residual}
To preserve the mathematical analogue to the strong variational multiscale formulation, we again choose a single bubble function for interpolating the fine-scale variables such that equation \eqref{eq:DefineBubbles} holds.  As usual, $\bar{\boldsymbol{v}}$ and its weighting function $\bar{\boldsymbol{w}}$ may be expressed in terms of the nodal values $\hat{\bar{\boldsymbol{v}}}$ and $\hat{\bar{\boldsymbol{w}}}$ as
\begin{align}
  \label{Eqn:SNS_Coarse_Interpolation}
\bar{\boldsymbol{v}} = \hat{\bar{\boldsymbol{v}}}^T \boldsymbol{N}^T \; ; \; \bar{\boldsymbol{w}} = \hat{\bar{\boldsymbol{w}}}^T \boldsymbol{N}^T
\end{align}
where $\boldsymbol{N}$ is a row vector of shape functions for each node.
Substituting equations and \eqref{eq:DefineBubbles} and \eqref{Eqn:SNS_Coarse_Interpolation} into equations \eqref{Eqn:SNS_Scalar_Resid_Coarse}-- \eqref{Eqn:SNS_Scalar_Resid_Fine} and noting the arbitrariness of $\hat{\bar{\boldsymbol{w}}}$ and $\boldsymbol{\gamma}$, we can construct vector residuals, $\boldsymbol{R}$, that are the sum contributions of the vector residuals at the element level given as
\begin{align}
  \label{Eqn:SNS_Coarse_Resid_Vec}
  \boldsymbol{R}^e_c(\bar{\boldsymbol{v}};\boldsymbol{v}',p) & := 2\nu \int_{\Omega_e}  \bar{\boldsymbol{B}}^T\mathrm{vec}[\nabla \bar{\boldsymbol{v}} + \nabla \boldsymbol{v}'] \; \mathrm{d} \Omega  - \int_{\Omega_e} \mathrm{vec}[\bar{\boldsymbol{G}}]p \; \mathrm{d} \Omega - \int_{\Omega_e} (\boldsymbol{N}^T\odot \boldsymbol{I})\boldsymbol{b} \; \mathrm{d} \Omega \\
  \label{Eqn:SNS_Coarse_Resid_Vec_2}
  \boldsymbol{R}^e_p(\bar{\boldsymbol{v}};\boldsymbol{v}') & := - \int_{\Omega_e} \boldsymbol{N}^T \nabla \cdot (\bar{\boldsymbol{v}} + \boldsymbol{v}') \; \mathrm{d} \Omega \\
  \label{Eqn:SNS_Fine_Resid_Vec}
  \boldsymbol{R}^e_f(\bar{\boldsymbol{v}};\boldsymbol{v}',p) & := 2\nu \int_{\Omega_e}  {\boldsymbol{B}'}^T
\mathrm{vec}[\nabla \bar{\boldsymbol{v}} + \nabla \boldsymbol{v}'] \; \mathrm{d} \Omega - \int_{\Omega_e}  {\boldsymbol{B}'}^T\mathrm{vec}[\boldsymbol{I}]p \; \mathrm{d} \Omega - \int_{\Omega_e} ( b^e \odot \boldsymbol{I})\boldsymbol{b} \; \mathrm{d} \Omega
  %
  %
  %
%
  %
  %
\end{align}
To write more compactly, we have made the substitutions
\begin{align}
\bar{\boldsymbol{B}} = \bar{\boldsymbol{G}} \odot \boldsymbol{I}&; \quad  \bar{\boldsymbol{G}} := \boldsymbol{J}^{-T}\boldsymbol{DN}^T \nonumber \\
\boldsymbol{B}' = \boldsymbol{g} \odot \boldsymbol{I}&; \quad  \boldsymbol{g} := \nabla_{\boldsymbol{x}} b^e
\end{align}
where $\boldsymbol{DN}$ represents a matrix of the first derivatives of the element shape functions,  $\boldsymbol{J}$ the element jacobian matrix, $\mathrm{vec}[\cdot]$ is an operation that represents a matrix with a vector, and $\odot$ is the Kronecker product \cite{Graham} (see Appendix). 
\subsection{Stiffness matrix}
Moving all applied force terms in $\boldsymbol{R}$ to the right hand side, we can write equations \eqref{Eqn:SNS_Coarse_Resid_Vec}--\eqref{Eqn:SNS_Fine_Resid_Vec} in matrix form as
\begin{align}
\label{Eqn:SNS_Matrix_Form_Resid}
  \boldsymbol{K} \left[ \begin{array}{c}  \bar{\boldsymbol{v}} \\
  p \\
  \boldsymbol{v}'  \end{array} \right] = \left[ \begin{array}{c}  \boldsymbol{f}_{c} \\
\boldsymbol{f}_{p} \\
\boldsymbol{f}_{f}  \end{array} \right]
\end{align}
where $\boldsymbol{f}$ represents the sum of the element contributions to the applied forces, defined as
\begin{align}
\boldsymbol{f}_{c}^e = \int_{\Omega_e} (\boldsymbol{N} \odot \boldsymbol{I}) \boldsymbol{b} \; \mathrm{d} \Omega ; \quad \boldsymbol{f}_{p}^e = 0; \quad \boldsymbol{f}_{f}^e = \int_{\Omega_e} ( b^e \odot \boldsymbol{I})\boldsymbol{b} \; \mathrm{d} \Omega
\end{align}
The global stiffness matrix, $\boldsymbol{K}$, before static condensation, has the form
\begin{align}
\label{eq:LinearStiffness}
\boldsymbol{K} = \left[ \begin{array}{ccc} \boldsymbol{K}_{cc} & \boldsymbol{K}_{cp} & \boldsymbol{K}_{cf} \\
 \boldsymbol{K}_{pc}  & \boldsymbol{K}_{pp}  & \boldsymbol{K}_{pf}  \\
 \boldsymbol{K}_{fc}  & \boldsymbol{K}_{fp}  & \boldsymbol{K}_{ff}   \end{array} \right] 
\end{align}
where the element contributions are computed as follows
\begin{align}
  \boldsymbol{K}_{cc}^e  &=  2\nu \int_{\Omega_e}  \bar{\boldsymbol{B}}^T \bar{\boldsymbol{B}} \; \mathrm{d} \Omega ; \quad
 \boldsymbol{K}_{cp}^e  =  - \int_{\Omega_e}  \mathrm{vec}[\bar{\boldsymbol{G}}]\boldsymbol{N} \; \mathrm{d} \Omega ; \quad
\boldsymbol{K}_{cf}^e  =  2\nu \int_{\Omega_e}  \bar{\boldsymbol{B}}^T \boldsymbol{B}'  \; \mathrm{d} \Omega \nonumber \\
 \boldsymbol{K}_{pc}^e  &=  - \int_{\Omega_e} \boldsymbol{N}^T  \mathrm{vec}[\bar{\boldsymbol{G}}]^T \; \mathrm{d} \Omega ; \quad
 \boldsymbol{K}_{pp}^e   =  0 ; \quad
\boldsymbol{K}_{pf}^e  =  - \int_{\Omega_e} \boldsymbol{N}^T {\boldsymbol{g}}^T  \; \mathrm{d} \Omega \nonumber \\
  \boldsymbol{K}_{fc}^e  &=  2\nu \int_{\Omega_e}  {\boldsymbol{B}'}^T\bar{\boldsymbol{B}}
 \; \mathrm{d} \Omega ; \quad
 \boldsymbol{K}_{fp}^e  =   - \int_{\Omega_e} {\boldsymbol{g}}^T \boldsymbol{N}\; \mathrm{d} \Omega ; \quad
 \boldsymbol{K}_{ff}^e  =  2\nu \int_{\Omega_e}  {\boldsymbol{B}'}^T\boldsymbol{B}'\; \mathrm{d} \Omega 
\end{align}

Using block Gauss elimination on equation \eqref{Eqn:SNS_Matrix_Form_Resid}, the fine-scale components can be condensed from the stiffness matrix.  The resulting matrix equation can be written as
\begin{align}
\label{Eqn:SNS_Condensed_Resid}
  \widetilde{\boldsymbol{K}} \left[ \begin{array}{c}  \bar{\boldsymbol{v}} \\
  p \end{array} \right] = \left[ \begin{array}{c}  \widetilde{\boldsymbol{f}}_{c} \\
\widetilde{\boldsymbol{f}}_{p} \end{array} \right]
\end{align}
The global stiffness matrix has the form
\begin{align}
\label{eq:CondensedLinearStiffness}
\widetilde{\boldsymbol{K}} = \left[ \begin{array}{cc} \widetilde{\boldsymbol{K}}_{cc} & \widetilde{\boldsymbol{K}}_{cp}  \\
 \widetilde{\boldsymbol{K}}_{pc}  & \widetilde{\boldsymbol{K}}_{pp} \end{array} \right] 
\end{align}
where we have augmented the coarse-scale components with the fine-scale components at the element level as follows
\begin{align}
\label{eq:LinearStiffness_Defined}
\widetilde{\boldsymbol{K}}_{cc}^e &= \boldsymbol{K}_{cc}^e - \boldsymbol{K}_{cf}^e(\boldsymbol{K}_{ff}^e)^{-1}\boldsymbol{K}_{fc}^e \\
\widetilde{\boldsymbol{K}}_{cp}^e &= \boldsymbol{K}_{cp}^e - \boldsymbol{K}_{cf}^e(\boldsymbol{K}_{ff}^e)^{-1}\boldsymbol{K}_{fp}^e \\
\widetilde{\boldsymbol{K}}_{pc}^e &= \boldsymbol{K}_{pc}^e - \boldsymbol{K}_{pf}^e(\boldsymbol{K}_{ff}^e)^{-1}\boldsymbol{K}_{fc}^e \\
\widetilde{\boldsymbol{K}}_{pp}^e &=  - \boldsymbol{K}_{pf}^e(\boldsymbol{K}_{ff}^e)^{-1}\boldsymbol{K}_{fp}^e 
\end{align}
Similarly, the applied force vector has been augmented at the element level as
\begin{align}
\label{eq:LinearForces_Defined}
  \widetilde{\boldsymbol{f}}_{c}^e = \boldsymbol{f}_{c}^e - \boldsymbol{K}_{cf}^e(\boldsymbol{K}_{ff}^e)^{-1}\boldsymbol{f}_{f}^e ; \quad \widetilde{\boldsymbol{f}}_{p}^e = - \boldsymbol{K}_{pf}^e(\boldsymbol{K}_{ff}^e)^{-1}\boldsymbol{f}_{f}^e
\end{align}
%

After solving for the coarse scale variables from equation \eqref{Eqn:SNS_Condensed_Resid}, the fine-scale variables can be recovered with post processing if desired.

%
\section{NUMERICAL RESULTS}

In this section we contrast the performance of the enriched formulation with that of the weak and strong variational multiscale stabilized formulations for various test problems. 

\subsection{Constant velocity and pressure problem}
\begin{figure}[htb!]
	\centering
  \includegraphics[scale=0.5]{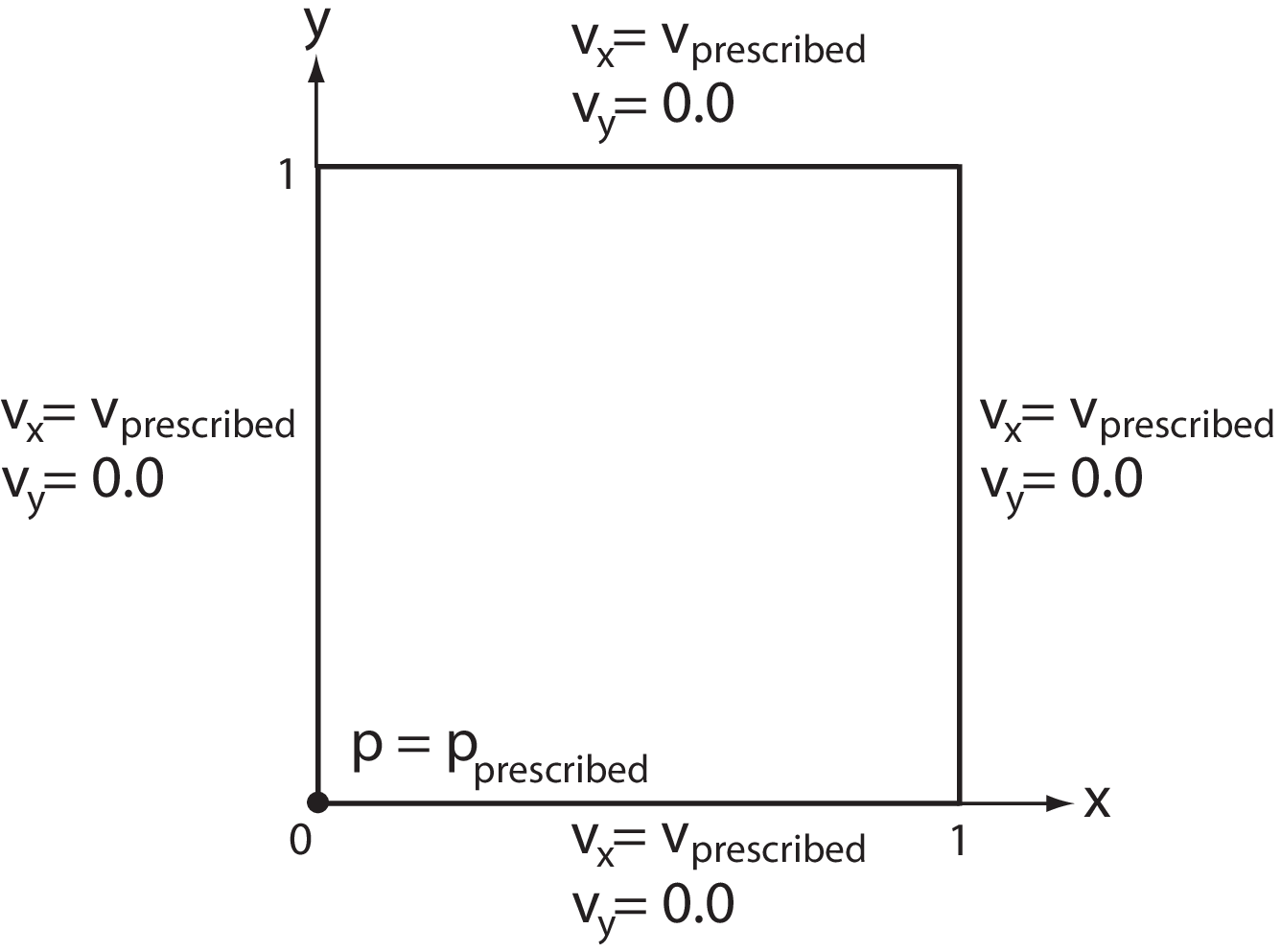}
	\caption{Constant velocity and pressure test problem: domain and boundary conditions.}
	\label{fig:ConstPressPatch}
\end{figure}

The constant velocity and pressure test problem represents an extremely simple physical state, yet even the most sophisticated formulation must capture it without oscillations. The solution to the constant velocity and pressure test problem is $\boldsymbol{v} = (10.0,0,0); p = 10.0$, which by inspection satisfies the governing equations (equations \eqref{Eqn:SNS_Equilibrium}--\eqref{Eqn:SNS_TractionBC}).  The boundary conditions are defined in Figure \ref{fig:ConstPressPatch} for the two-dimensional case.
\subsubsection{TET4 elements.}
\begin{figure}[htb!]
	\centering
  \subfigure[]{\includegraphics[scale=0.30]{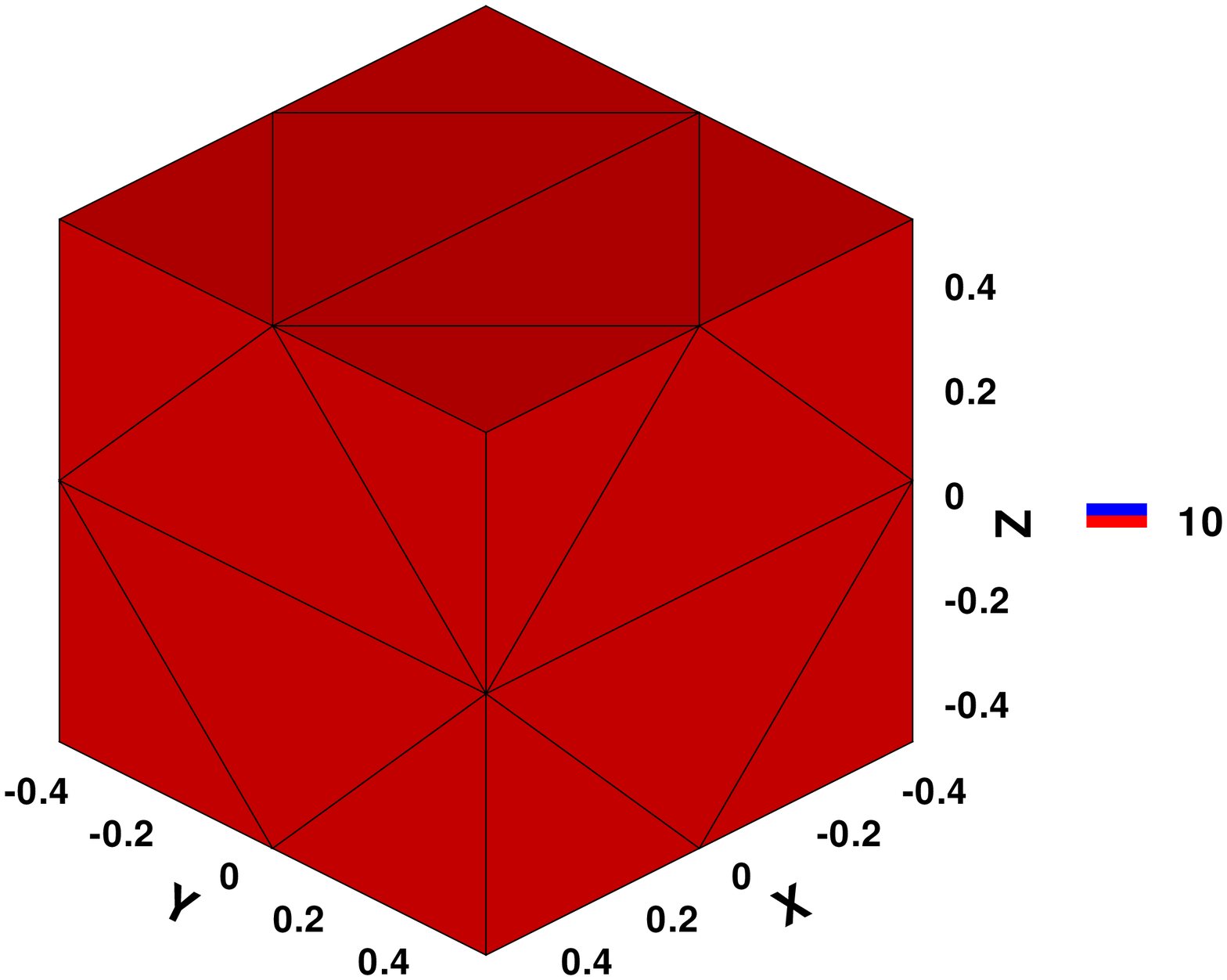}}
  \subfigure[]{\includegraphics[scale=0.30]{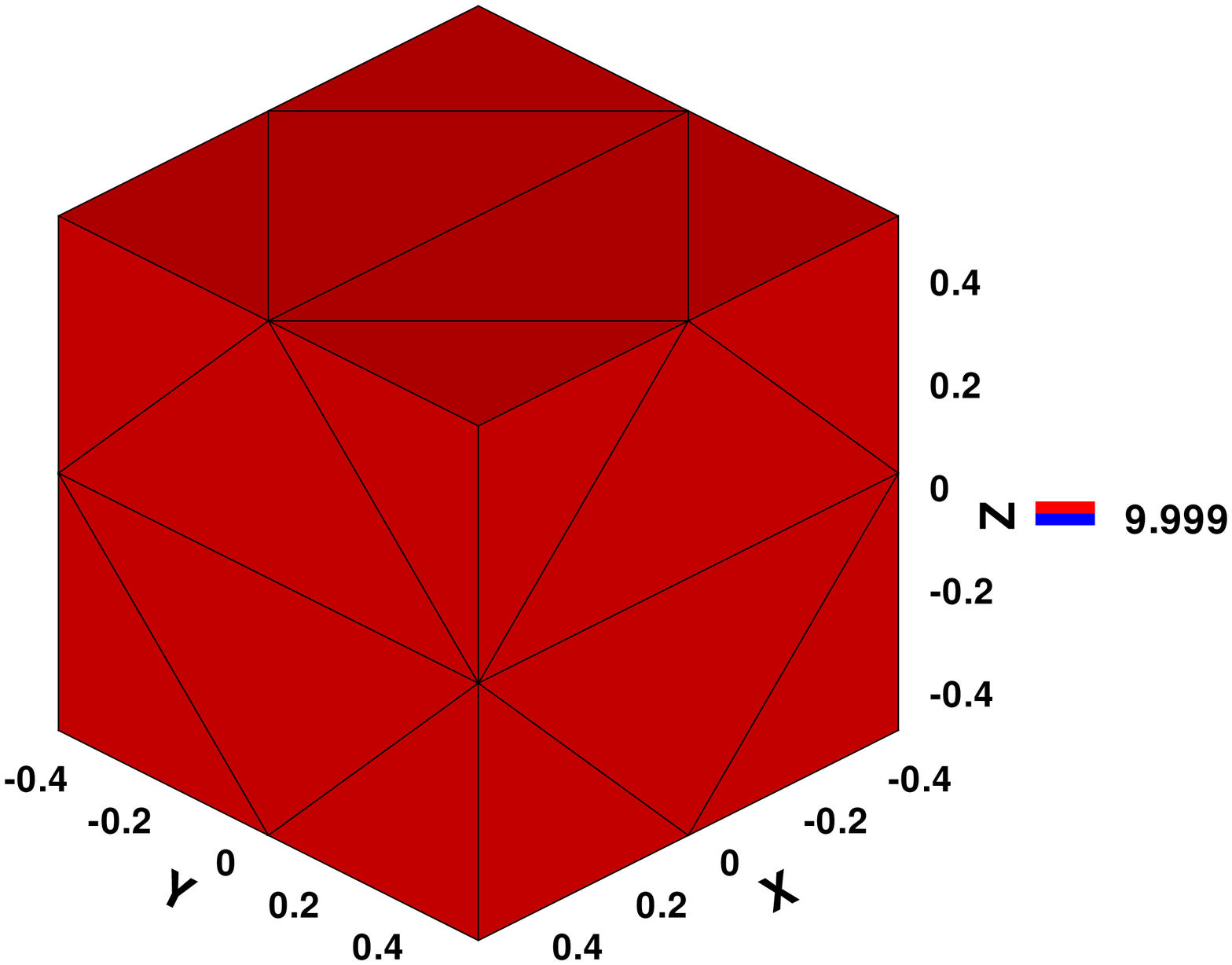}}
	\caption{Constant velocity and pressure test problem: (a) x-velocity and (b) pressure for 36 TET4 elements using the enriched formulation.}
	\label{fig:TET4VelPress}
\end{figure}

\begin{figure}[htb!]
	\centering
  \subfigure[]{\includegraphics[scale=0.30]{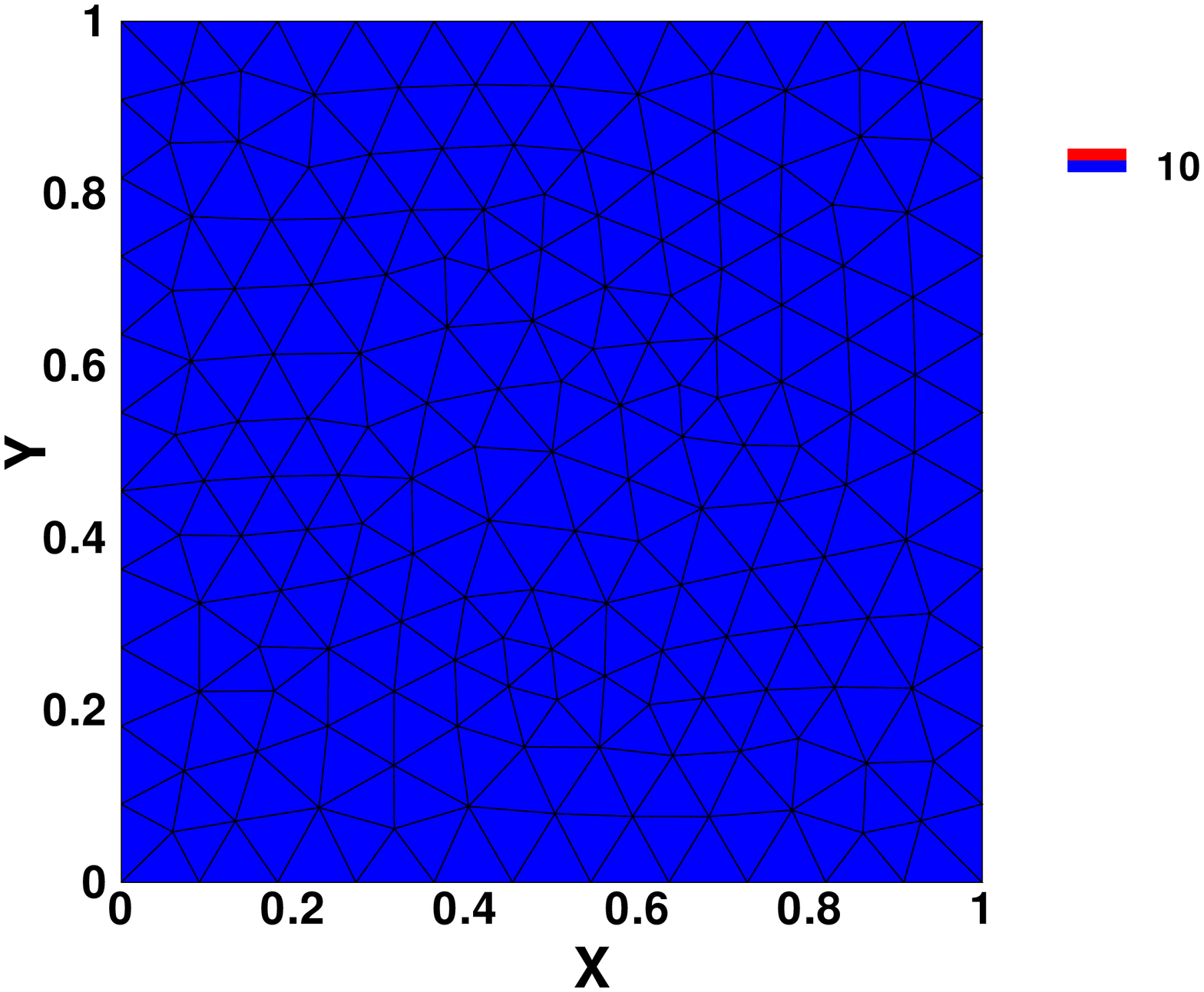}}
  \subfigure[]{\includegraphics[scale=0.30]{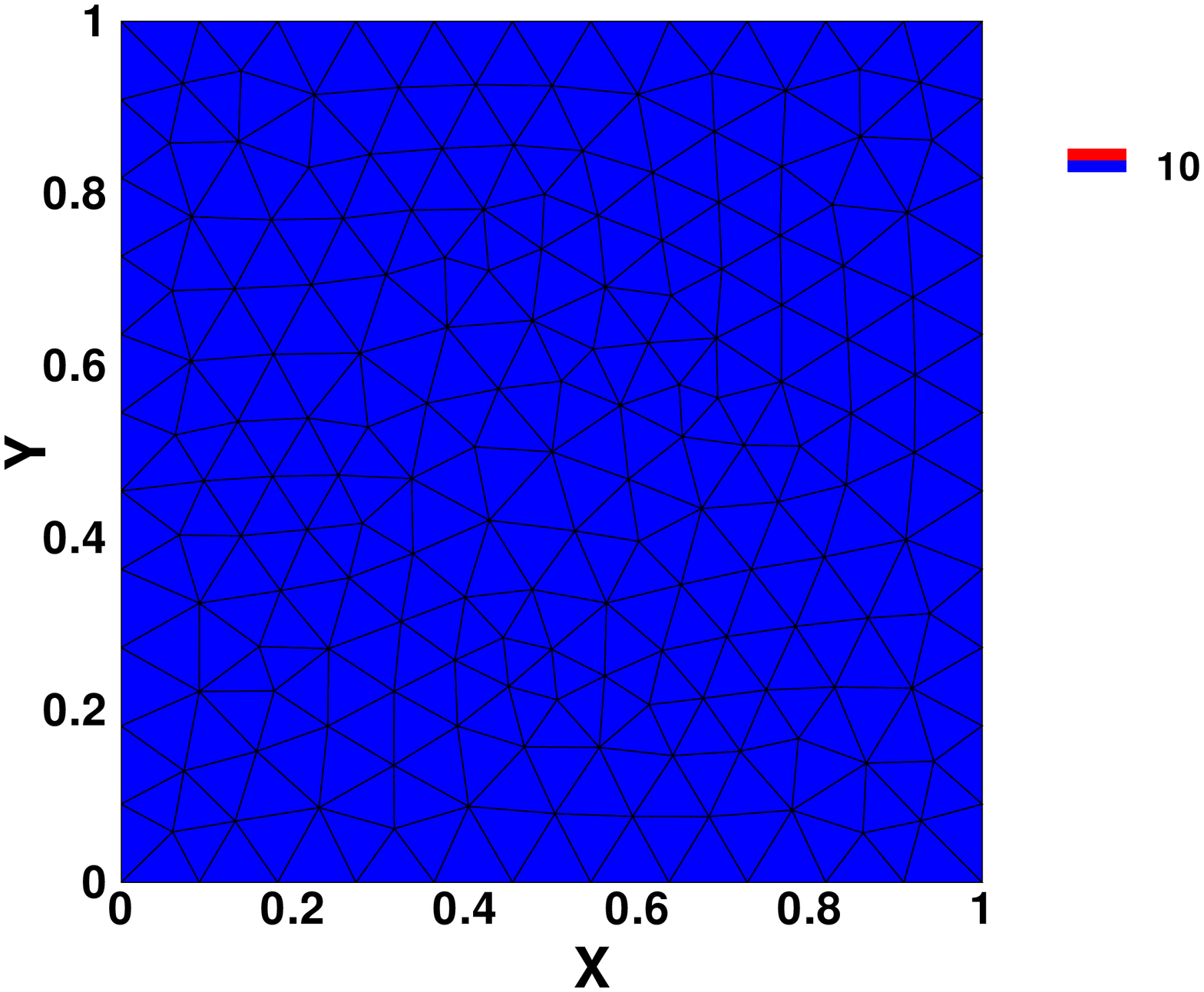}}
	\caption{Constant velocity and pressure test problem: (a) x-velocity and (b) pressure for 336 well-centered triangle elements using a standard Galerkin formulation.}
	\label{fig:WCTVelPress}
\end{figure}

As already mentioned in the introduction, the stability of the enriched method has been proved for triangular elements, but for the sake of completeness, we show that TET4 elements also perform well for the constant velocity and pressure problem.  The results are shown in Figure \ref{fig:TET4VelPress}.
\begin{remark}
As an aside, the authors would also like to point out that for a well-centered triangle (WCT) mesh (triangles with no interior angles greater than or equal to 90 degrees), even the standard Galerkin formulation produces no oscillations for the constant velocity and pressure problem.  The results are shown in Figure \ref{fig:WCTVelPress}.  A proof for the stability of such meshes is yet to appear in the literature.
\end{remark}
%
%
%

\subsubsection{Q4 elements.}
\begin{figure}[htb!]
	\centering
  \subfigure[]{\includegraphics[scale=0.30]{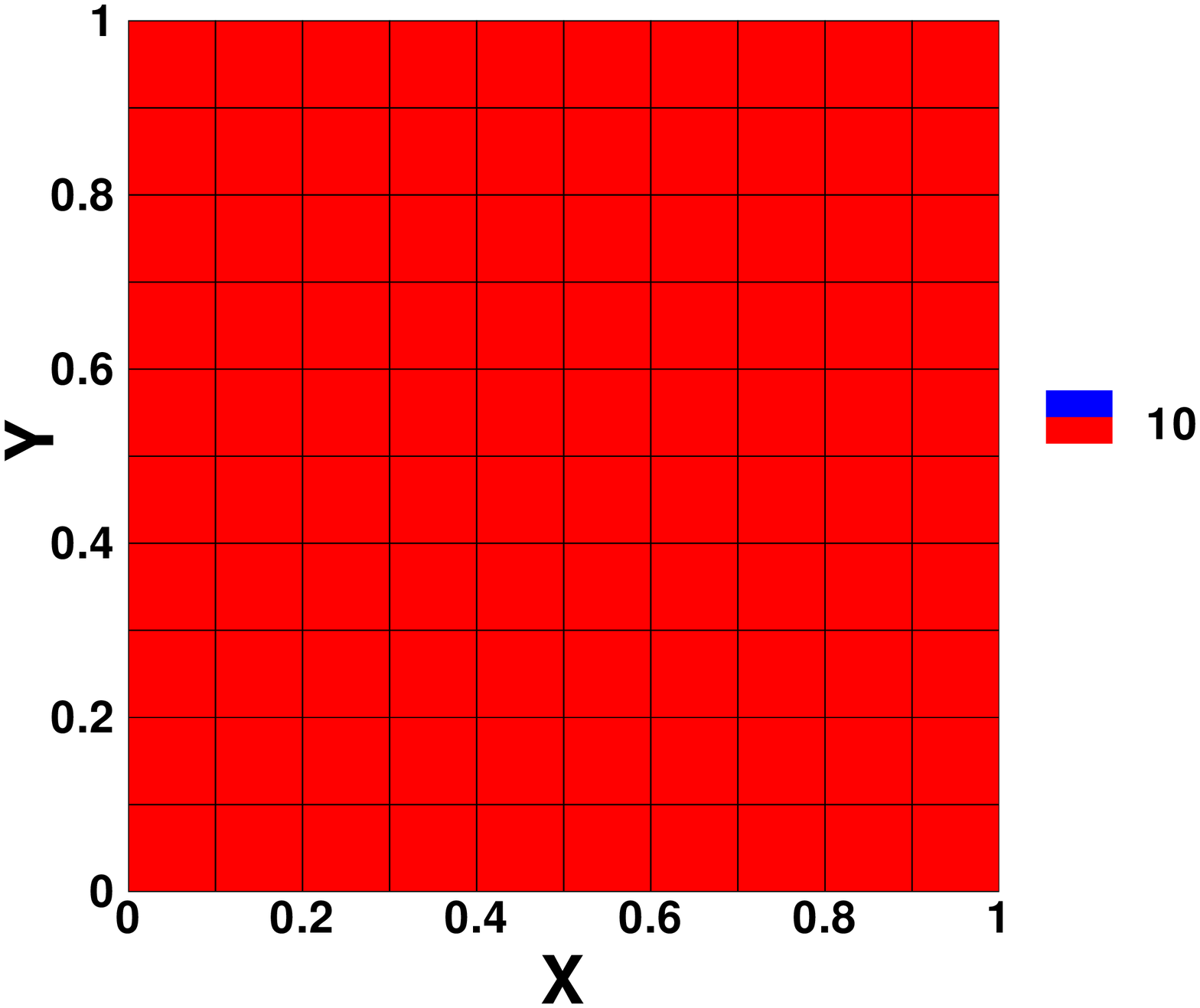}}
  \subfigure[]{\includegraphics[scale=0.30]{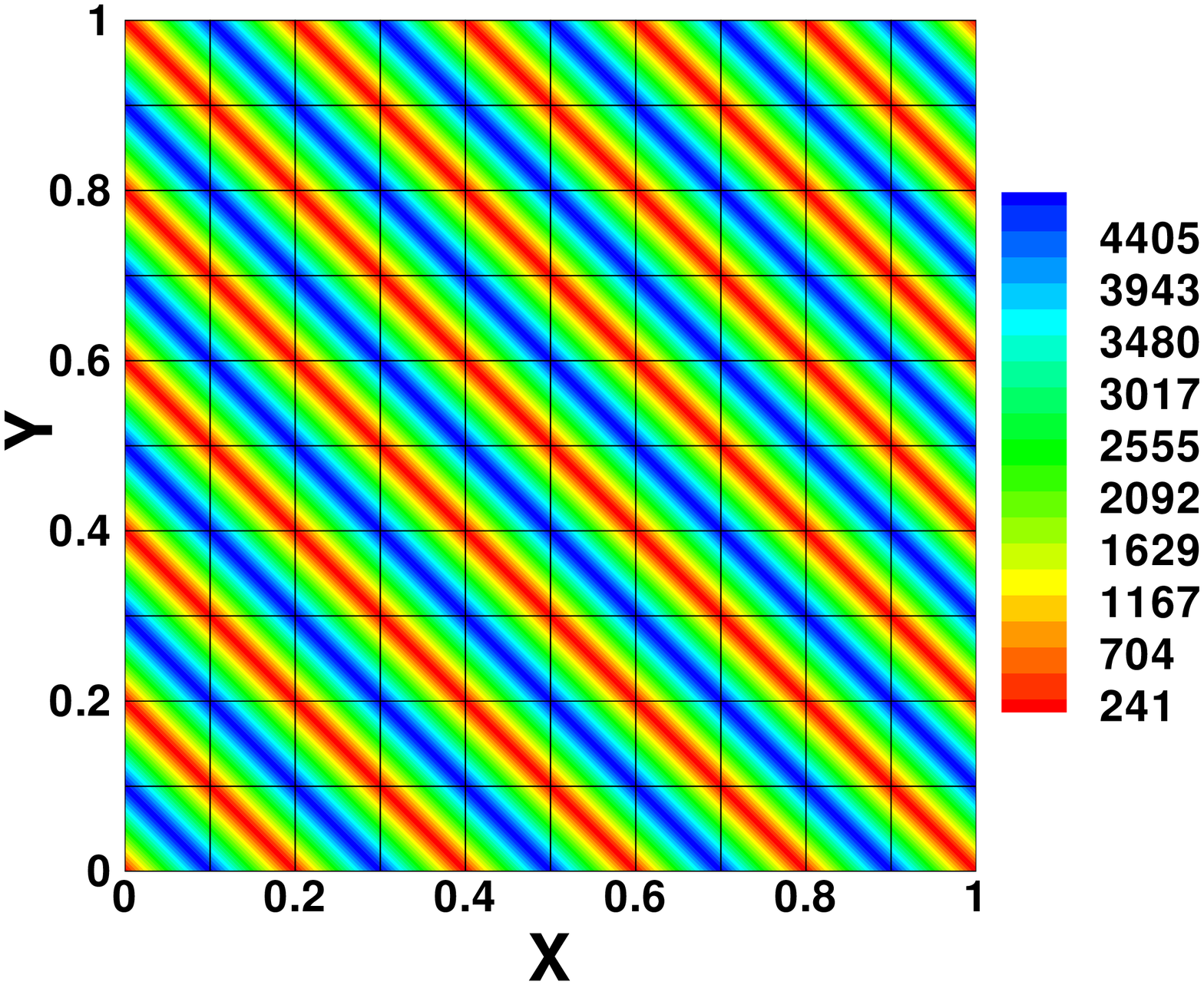}}
	\caption{Constant velocity and pressure test problem: (a) x-velocity and (b) pressure for 100 Q4 elements using the enriched formulation.}
	\label{fig:Kabab1}
\end{figure}

As pointed out in Remark \ref{remark:ResidRemark}, the statement that the coarse-scale residual is constant over the element domain in the limit of refinement is exactly true for T3 and TET4 elements for a constant body force, but in the case of Q4 and B8 elements, this statement is only an approximation.  Due to the introduction of this approximation, the enriched and stabilized formulations produce starkly contrasting results for Q4 elements when applied to the constant velocity and pressure problem.  Neither the weak or strong variational multiscale formulation shows any oscillations in the pressure or velocity, but as shown in Figure \ref{fig:Kabab1}, one can see that the enriched formulation shows severe pressure oscillations.
Brezzi and Pitk\"{a}ranta \cite{Mulder} proposed a stabilizing technique to remedy such spurious modes by circumventing the LBB condition.  To do so, one augments the enriched formulation with an added stability term $\epsilon(\nabla q;\nabla p)$ where $\epsilon \approx \mathcal{O}(h^2)$.  This resolves the ``missing'' $\boldsymbol{K}_{pp}$ term in the stiffness matrix before static condensation.  Performing this augmentation indeed weakly stabilizes the constant velocity and pressure problem, but this artificial term is not mathematically consistent.  A similar approach, the Pressure Stabilizing/Petrov-Galerkin (PSPG) method \cite{PSPG} circumvents the LBB condition, but preserves consistency by applying a perturbed weight function to all terms in the momentum equation.  Although the PSPG method avoids oscillations in the pressure, the stability parameter is usually defined in an ad hoc manner.

\begin{figure}[htb!]
	\centering
  \subfigure[]{\includegraphics[scale=0.35]{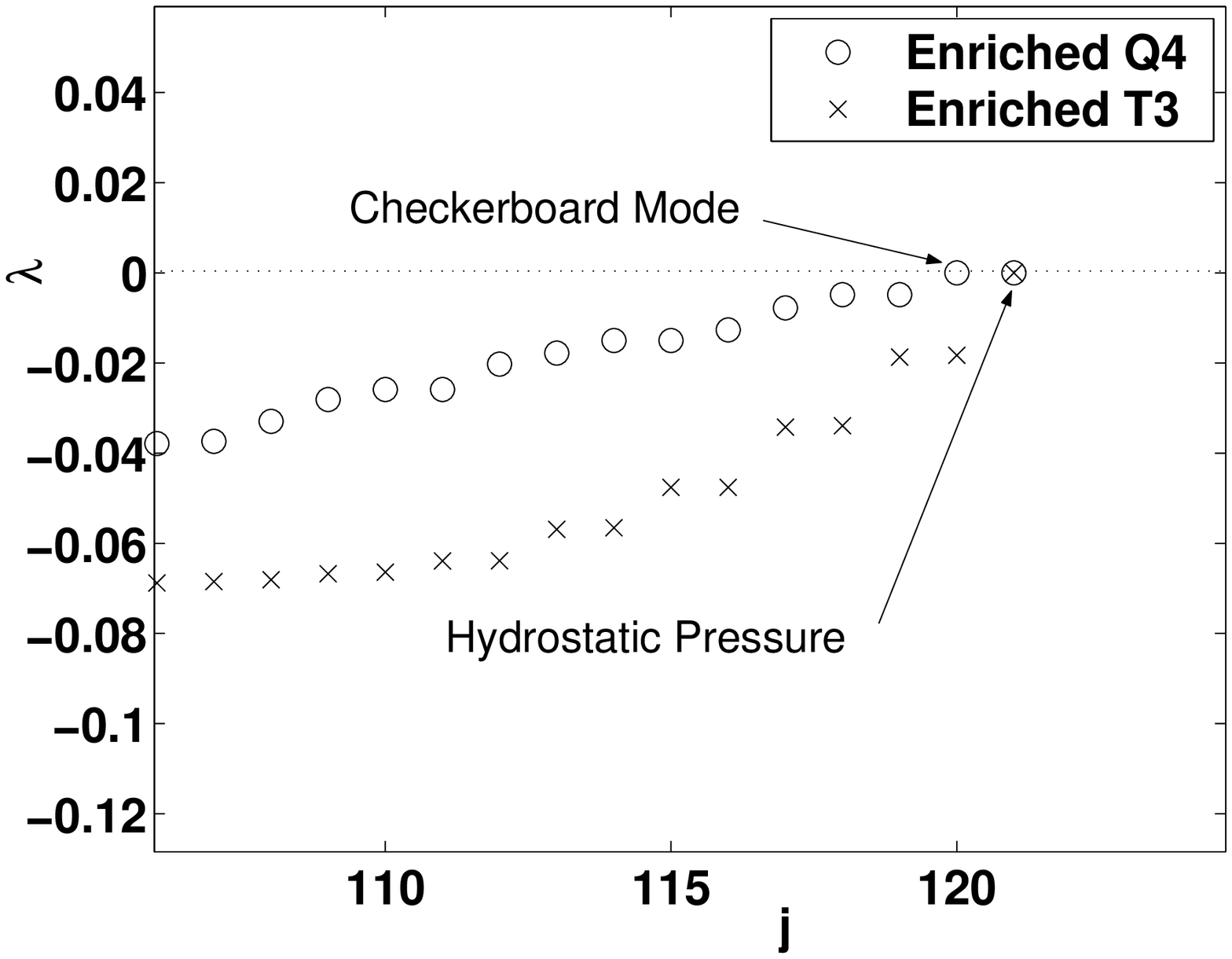}}
  \subfigure[]{\includegraphics[scale=0.35]{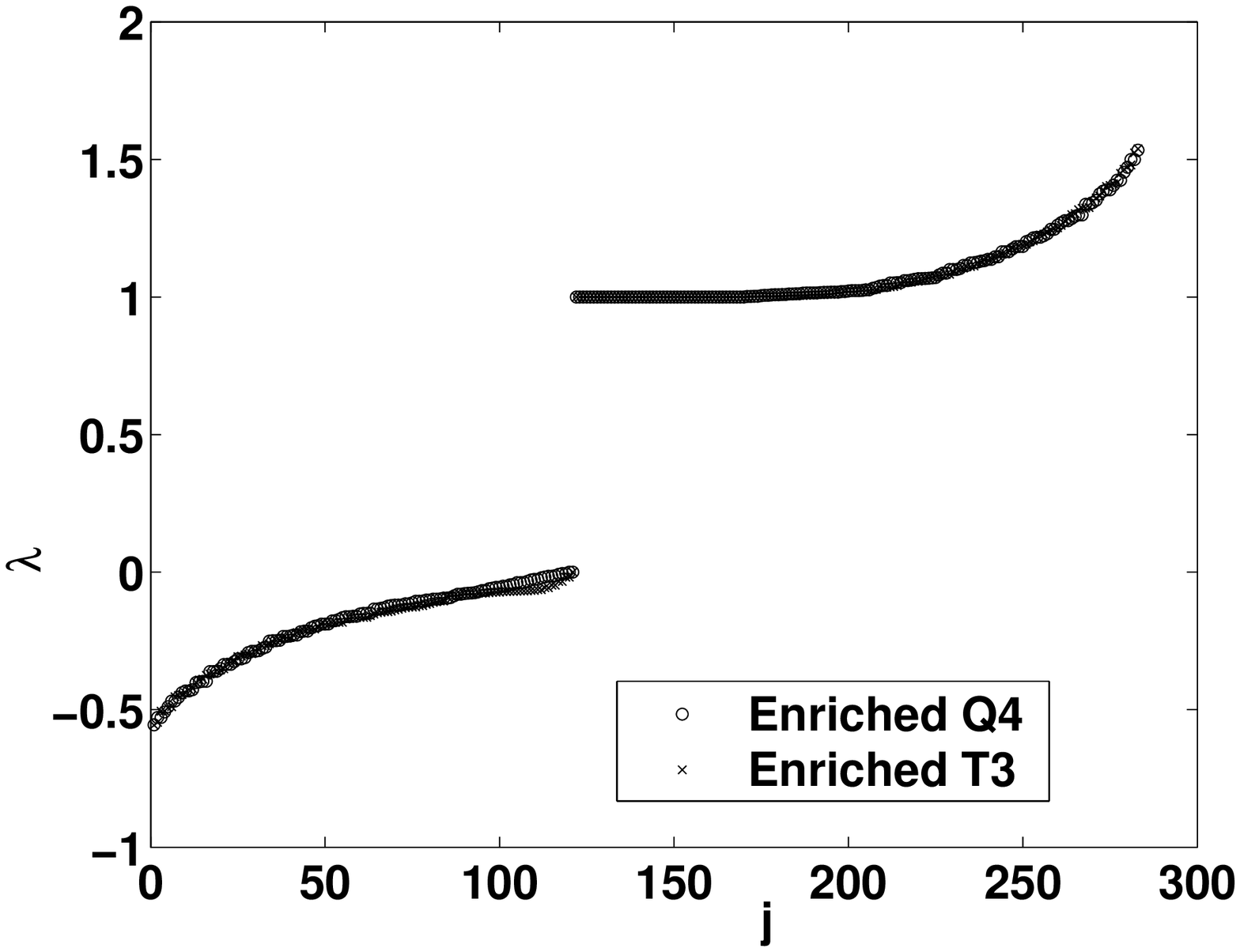}}
	\caption{Eigenvalues, $\boldsymbol{\lambda}$, associated with the discrete Stokes problem for enriched Q4 vs. enriched T3 (MINI) elements for $1/h$ = 10 (a) close-up of pure pressure modes (b) all eigenvalues shown.}
	\label{fig:Eigenvalue}
\end{figure}

In \cite{Malkus, Griffiths, Gresho, Sani}, the authors present an eigenvalue problem associated with the discrete LBB condition.  Analysis of the eigenvalue spectrum reveals certain oscillatory modes, for example the pure pressure modes for which the associated eigenvalues are zero.  The pure pressure modes consist of the hydrostatic mode and the checkerboard mode.  The hydrostatic mode can easily be removed by properly prescribing the pressure boundary conditions, but the checkerboard mode is related to linear dependence in the discretized system of equations. Using a unit square discretized with a grid of $n \times n$ enriched Q4 and T3 elements, we present the results from a similar eigenvalue analysis in Figure \ref{fig:Eigenvalue}.  The results show that bubble enrichment removes the checkerboard mode for the T3 (MINI) element, but that the checkerboard mode remains for the enriched Q4 element.  The presence of the checkerboard mode for the enriched Q4 element is consistent with the results shown in Figure \ref{fig:Kabab1}.
%
%
\subsubsection{B8 elements.}
\begin{figure}[htb!]
	\centering
  \subfigure[]{\includegraphics[scale=0.30]{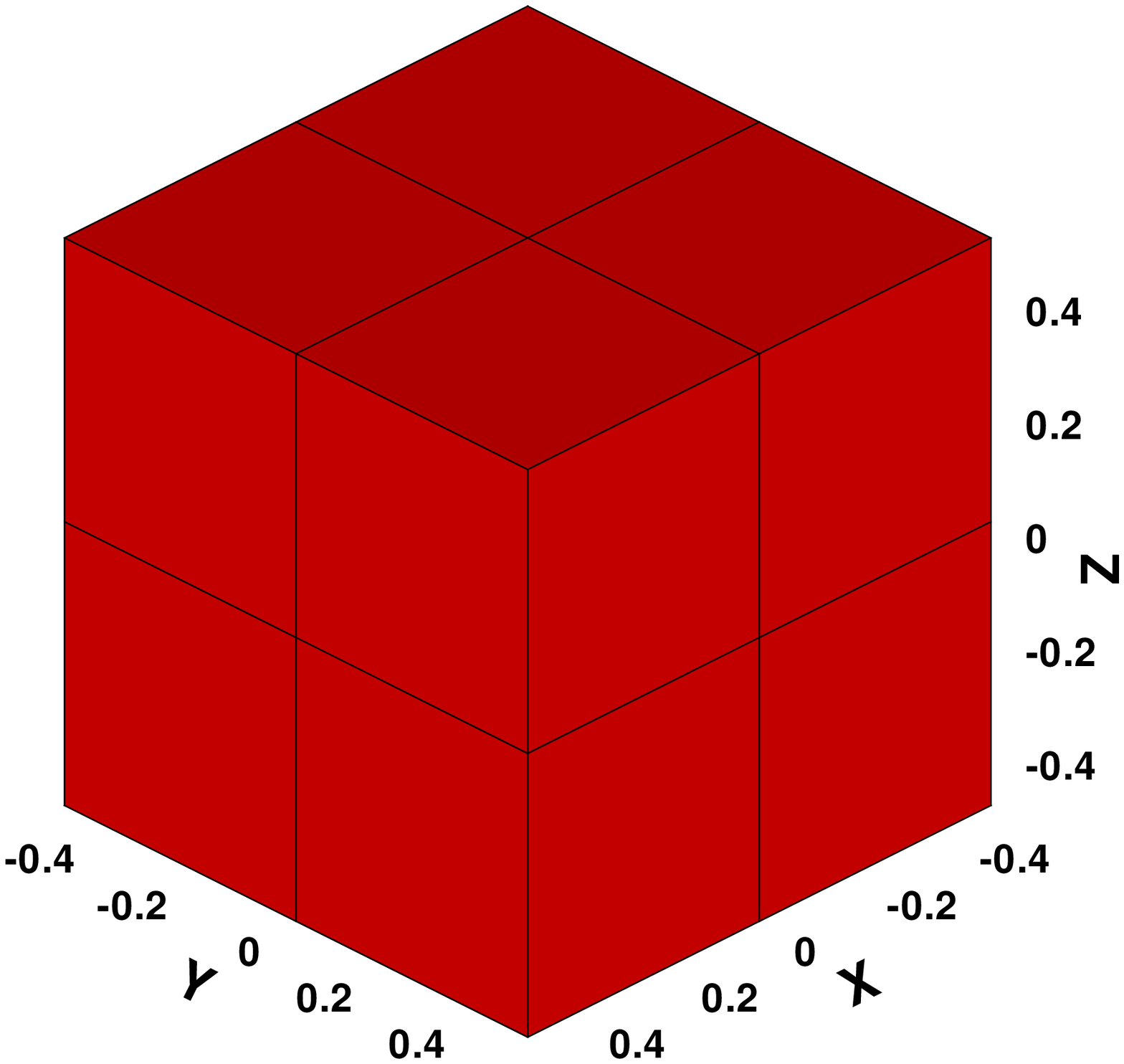}}
  \subfigure[]{\includegraphics[scale=0.30]{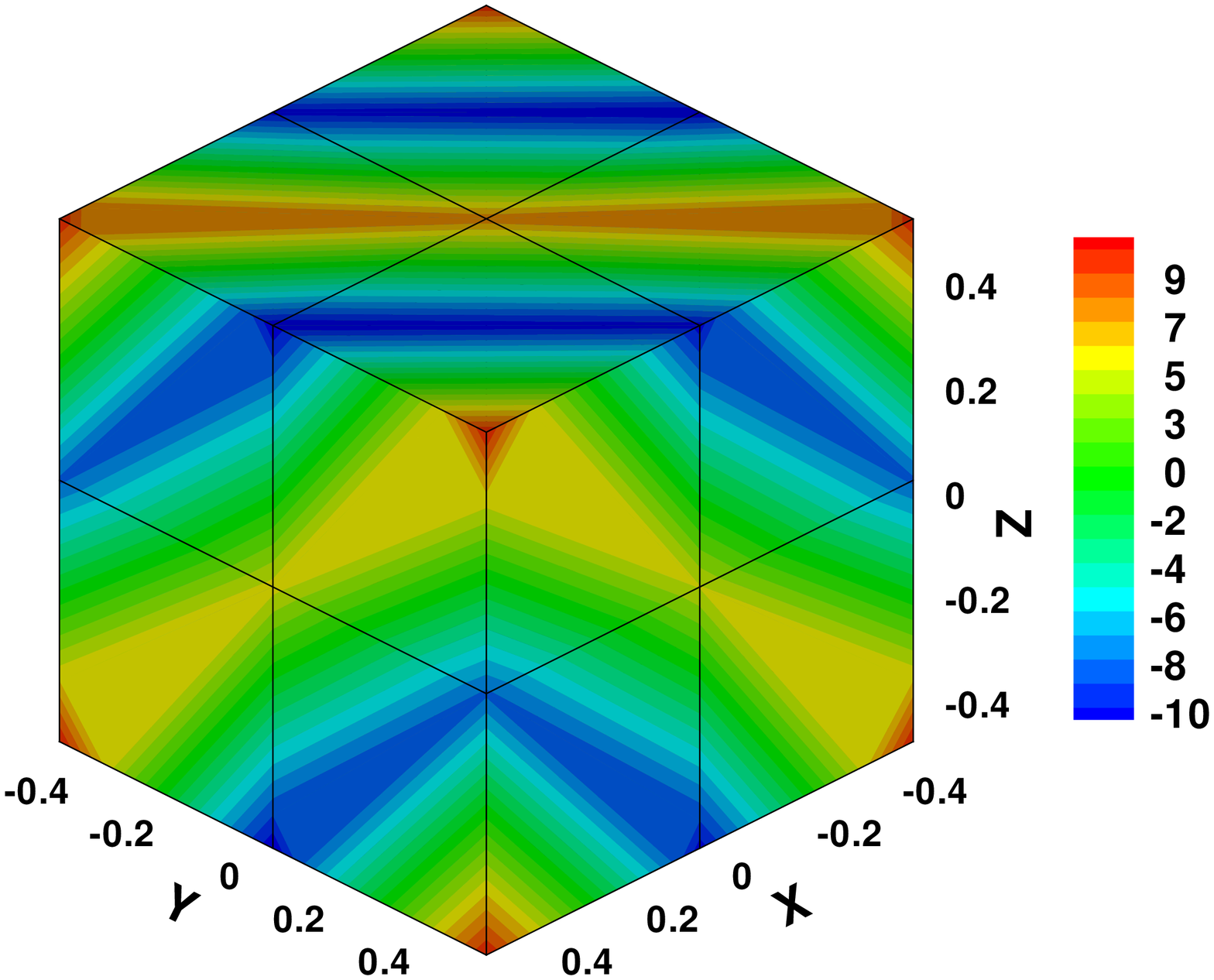}}
	\caption{Constant velocity and pressure test problem: (a) x-velocity and (b) pressure for 8 B8 elements using the enriched formulation.}
	\label{fig:ConstantFlowCoarse}
\end{figure}

\begin{figure}[htb!]
	\centering
  \subfigure[]{\includegraphics[scale=0.30]{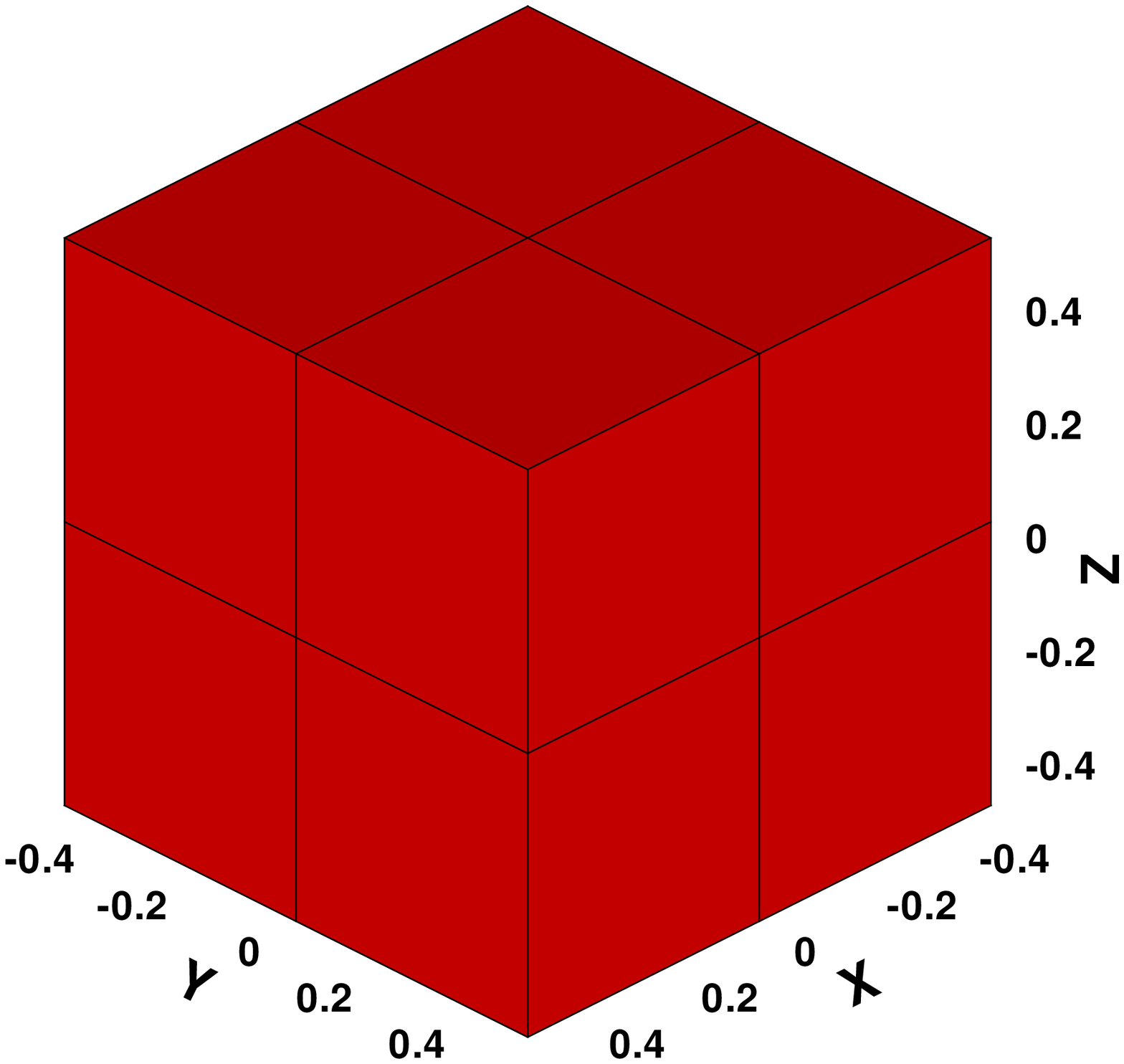}}
  \subfigure[]{\includegraphics[scale=0.30]{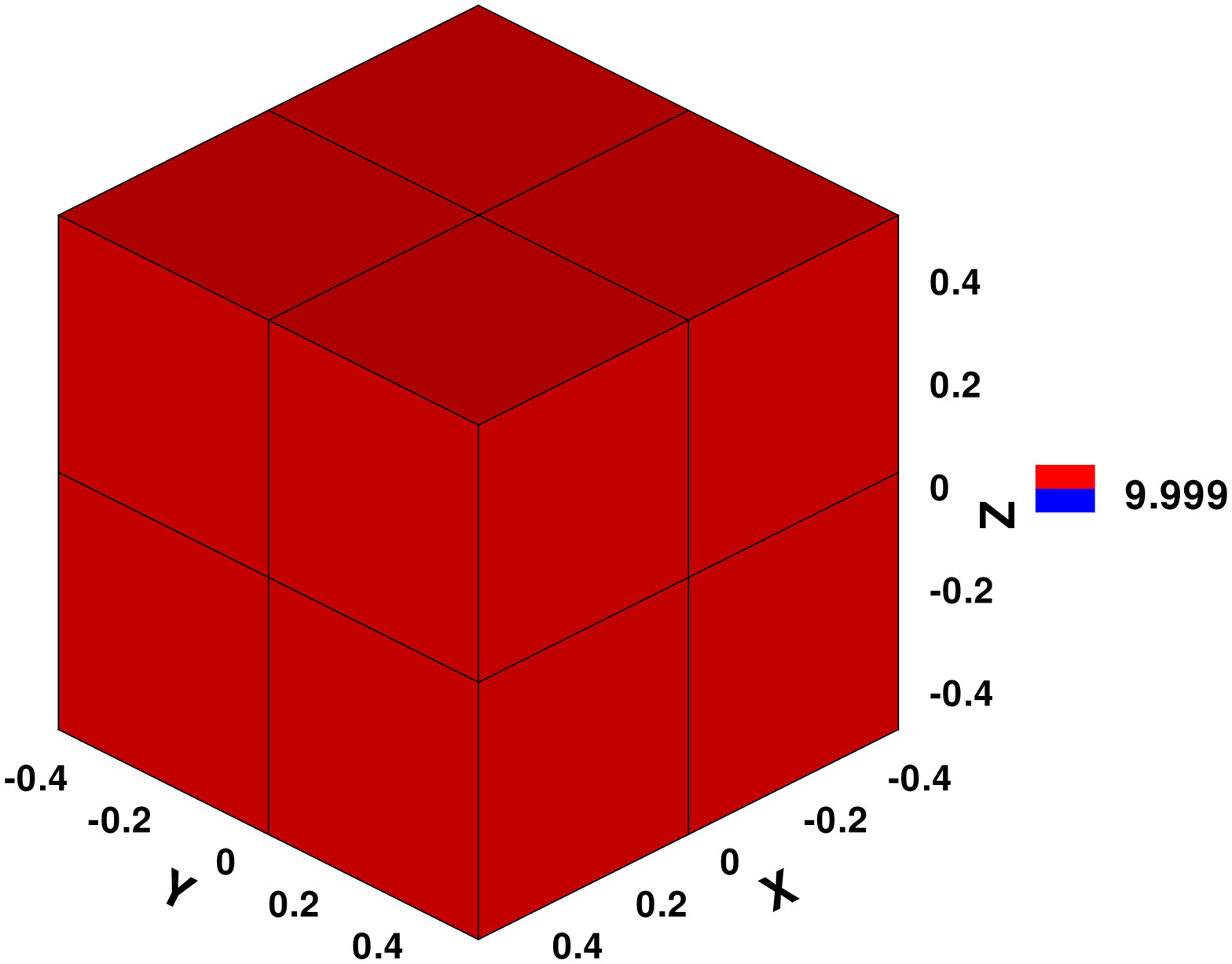}}
	\caption{Constant velocity and pressure test problem: (a) x-velocity and (b) pressure for 8 B8 elements using the weak variational multiscale formulation.}
	\label{fig:StabilizedCoarse}
\end{figure}

\begin{figure}[htb!]
	\centering
  \subfigure[]{\includegraphics[scale=0.30]{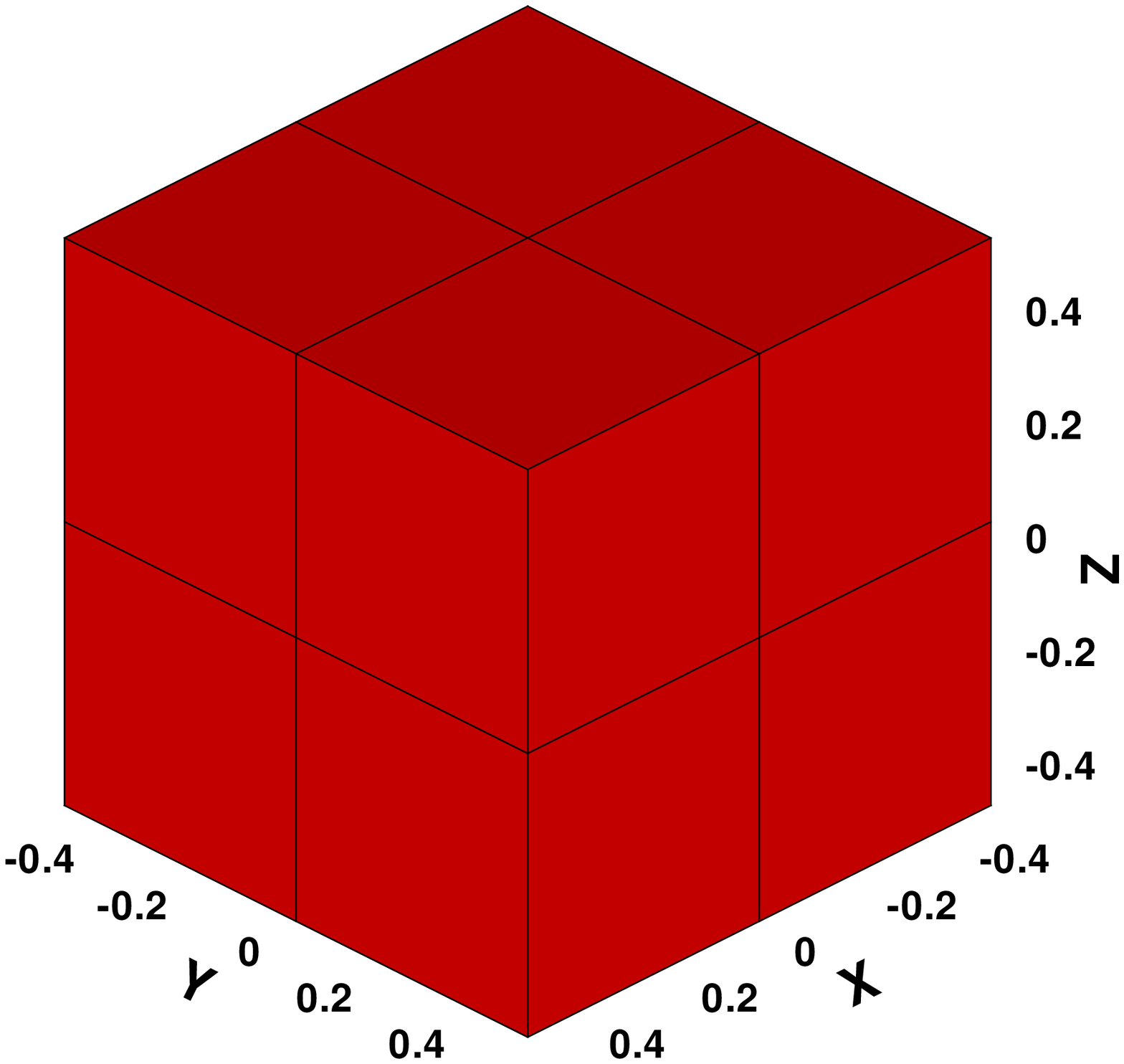}}
  \subfigure[]{\includegraphics[scale=0.30]{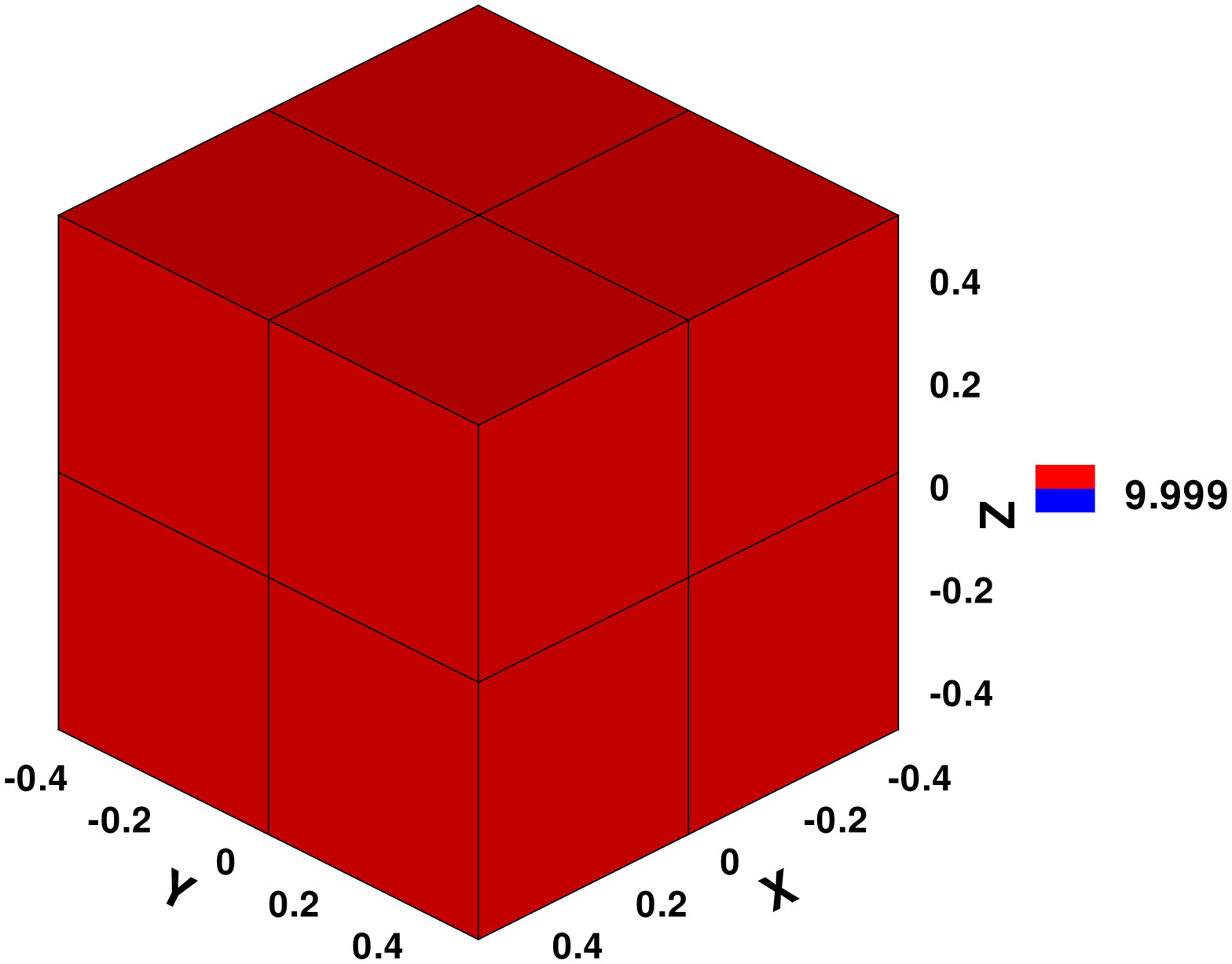}}
	\caption{Constant velocity and pressure test problem: (a) x-velocity and (b) pressure for 8 B8 elements using the strong variational multiscale formulation.}
	\label{fig:StabilizedCoarse2}
\end{figure}

\begin{figure}[htb!]
	\centering
  \subfigure[]{\includegraphics[scale=0.30]{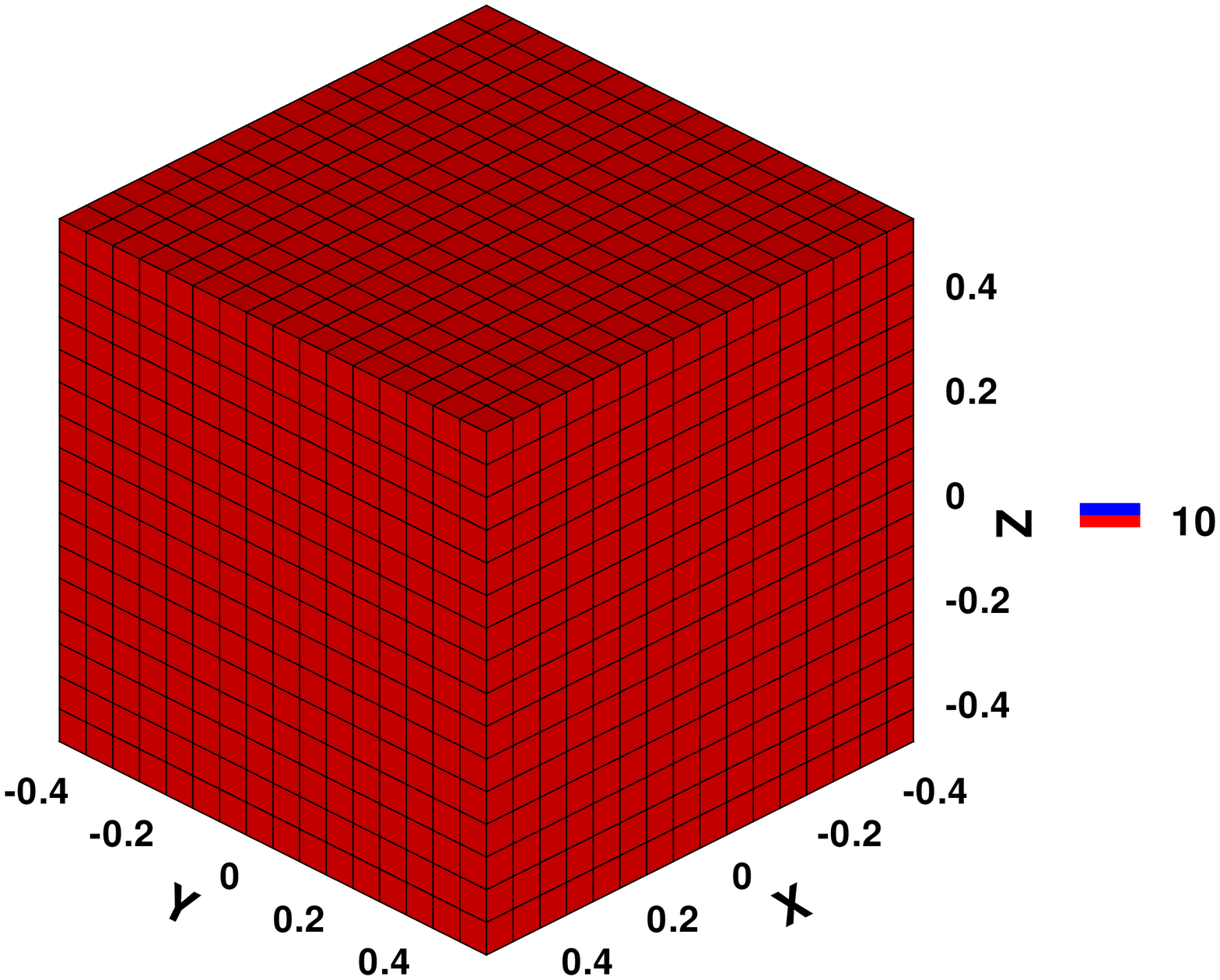}}
  \subfigure[]{\includegraphics[scale=0.30]{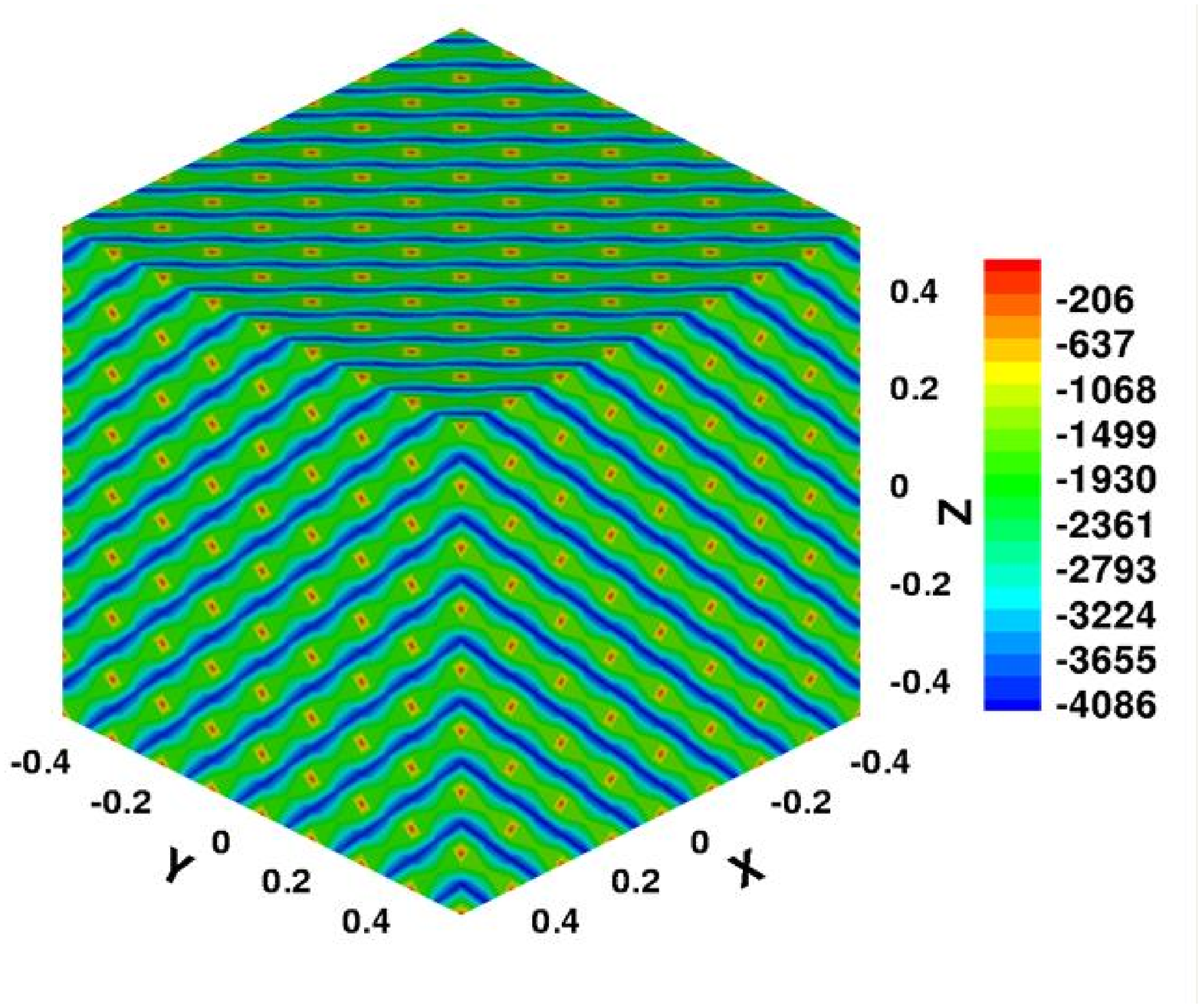}}
	\caption{Constant velocity and pressure test problem: (a) x-velocity and (b) pressure for 4096 B8 elements using the enriched formulation.}
	\label{fig:ConstantFlowRefined}
\end{figure}

Results similar to the two-dimensional case are obtained when extended to three-dimensions.  In particular, the B8 element also shows non-physical oscillations for this test problem that increase with mesh refinement.  Figure \ref{fig:ConstantFlowCoarse} shows the results of the three-dimensional test problem for a coarse mesh and Figure \ref{fig:ConstantFlowRefined} shows the results for a refined mesh.  Notice that no oscillations are present for the results obtained with the weak or strong variational multiscale  formulations as presented in Figures \ref{fig:StabilizedCoarse} and \ref{fig:StabilizedCoarse2}.
\subsection{Strong variational multiscale formulation}
To further verify the strong variational multiscale formulation of the Stokes problem, we present the results for some test problems along with a convergence analysis.  
\subsubsection{Lid--driven cavity.}
\begin{figure}[htb!]
	\centering
  \includegraphics[scale=0.6]{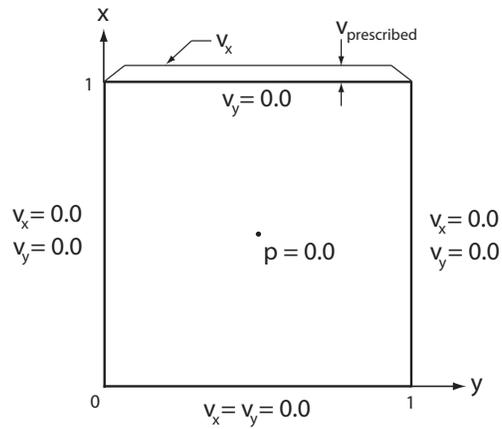}
	\caption{Lid--driven cavity: problem statement and boundary conditions.  The non-leaky cavity approach is used here which resolves the discontinuity at the upper two corners of the domain by assuming that the corners belong to the vertical walls.}
	\label{fig:Cavity}
\end{figure}

\begin{figure}[htb!]
	\begin{center}
  \subfigure[]{\includegraphics[scale=0.30]{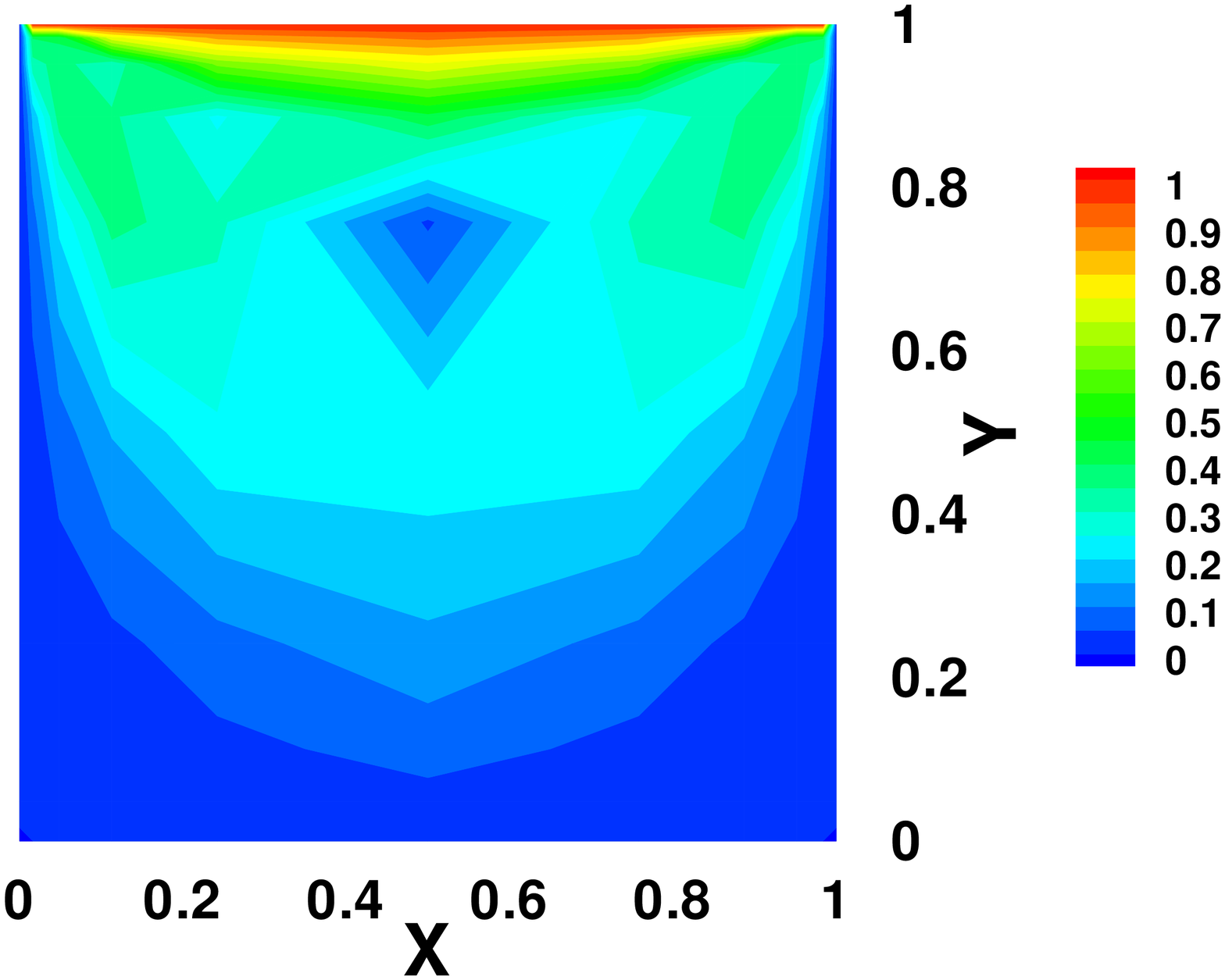}}
  \subfigure[]{\includegraphics[scale=0.30]{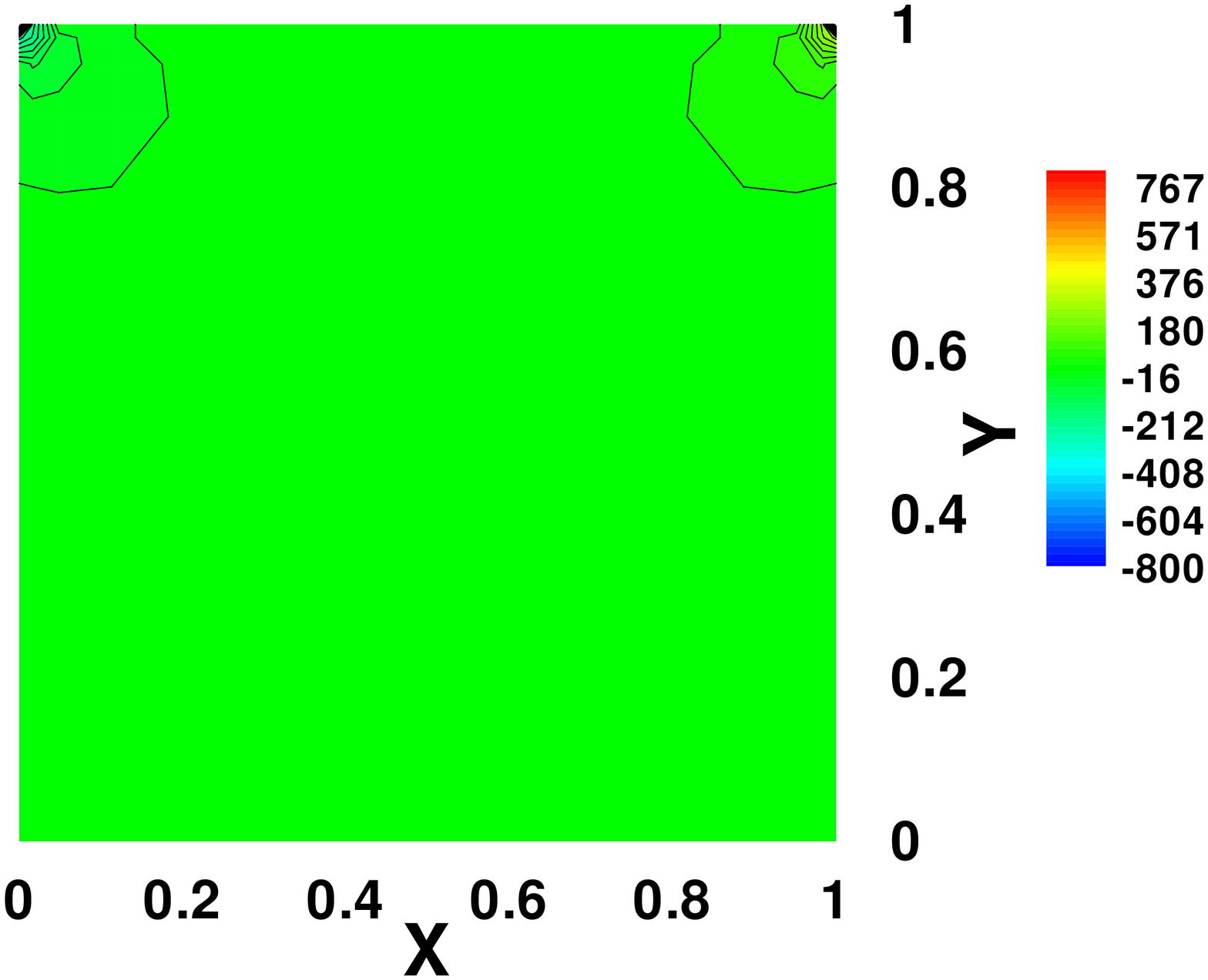}}
  \end{center}
	\caption{Lid--driven cavity test problem: (a) velocity magnitude and (b) pressure for 100 B8 elements.}
	\label{fig:CavityVelPress}
\end{figure}

The first problem evaluated is the well-known Lid-driven cavity problem.  A description of the domain, along with the boundary conditions is shown in Figure \ref{fig:Cavity}.
Contours of the velocity and pressure are shown in Figure \ref{fig:CavityVelPress}.  The results are in good accordance with other published results as shown in Table \ref{table:VortexLocation}.
\begin{table}[htb]
\caption{Position of the main cavity vortex}
\centering
\begin{tabular}{l c c }
\hline
  & element type & y-location  \\ [0.5ex]
\hline
Present simulation  & B8 & 0.753    \\
Donea and Huerta\cite{Donea} & Q2Q1 & 0.756 \\
\hline
\end{tabular}
\label{table:VortexLocation}
\end{table}
\subsubsection{Body force driven cavity.}
\begin{figure}[htb!]
	\begin{center}
  \subfigure[]{\includegraphics[scale=0.30]{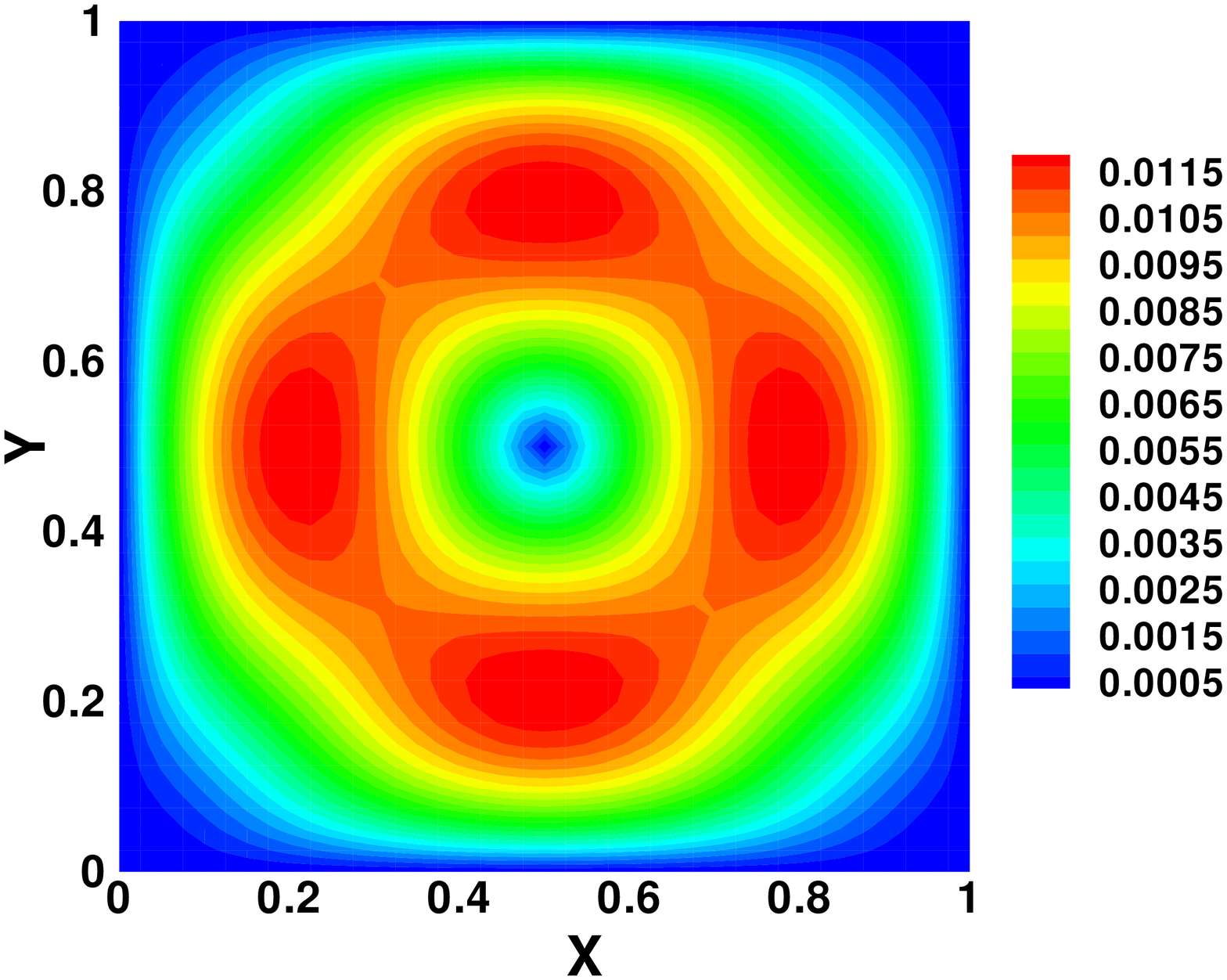}}
  \subfigure[]{\includegraphics[scale=0.30]{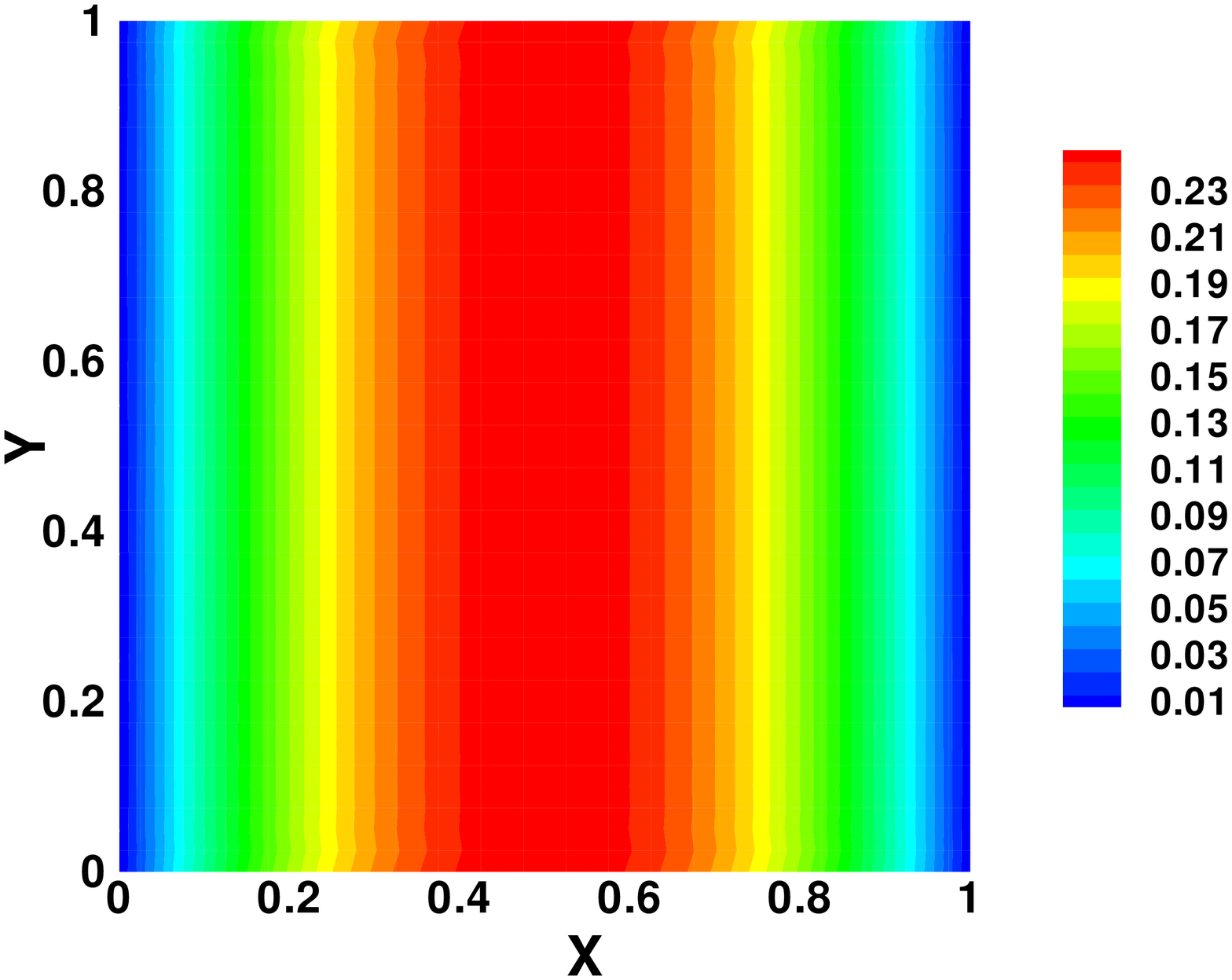}}
  \end{center}
	\caption{Body force driven cavity test problem: (a) velocity magnitude and (b) pressure computed with 1600 Q4 elements.}
	\label{fig:BodyForce}
\end{figure}

\begin{figure}[htb!]
	\centering
	\subfigure[]{\includegraphics[scale=0.35]{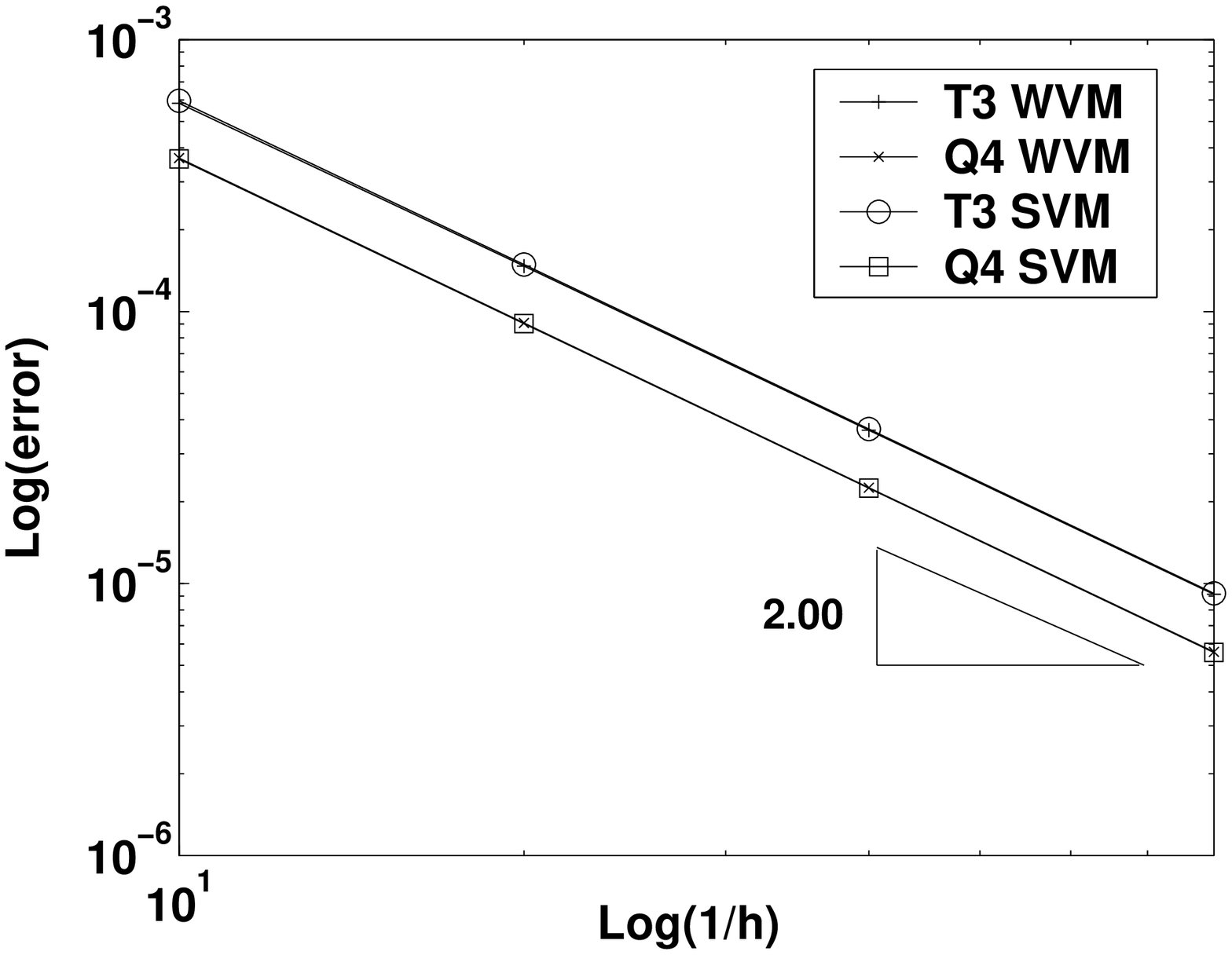}}
  \subfigure[]{\includegraphics[scale=0.35]{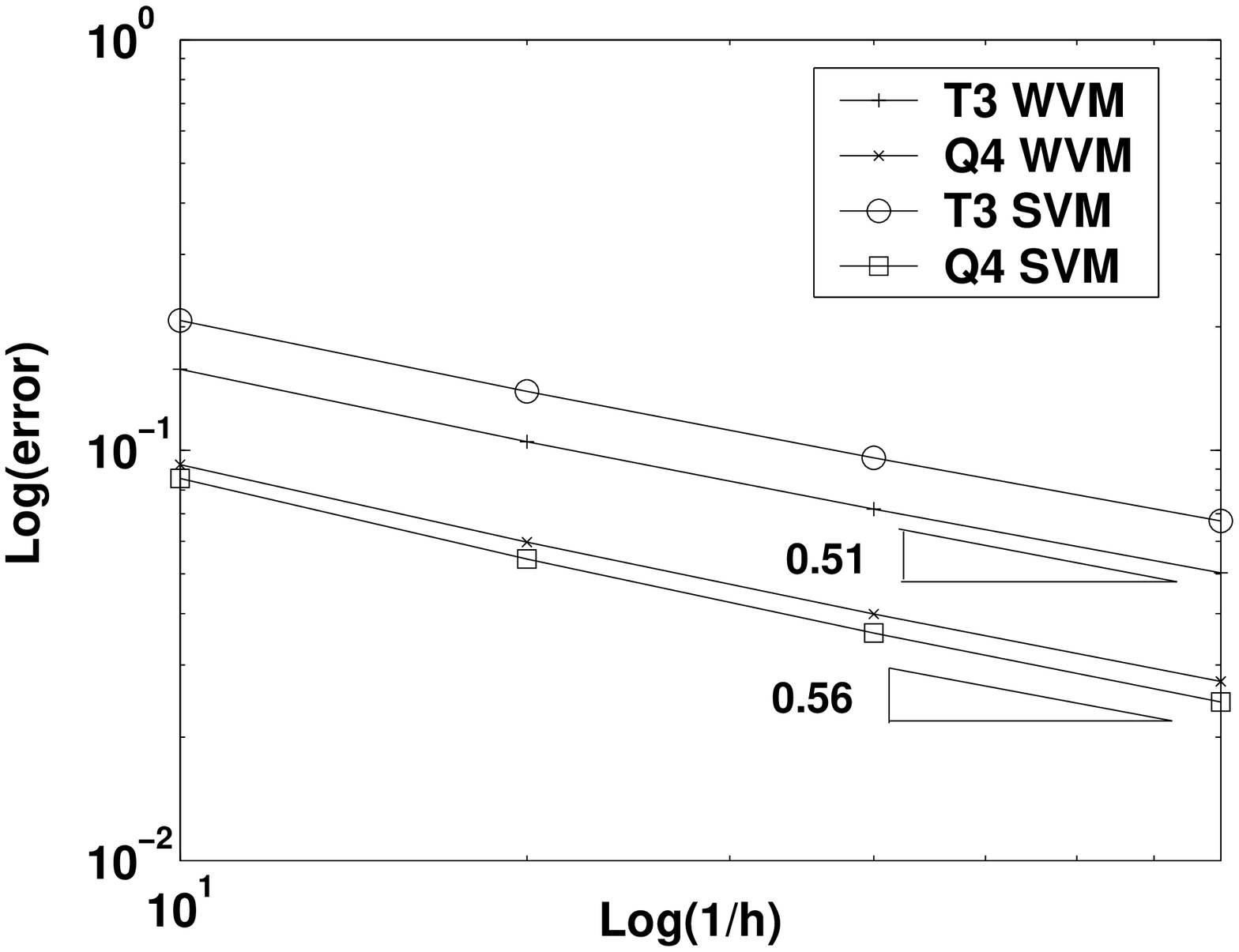}}
	\caption{Convergence rates on a uniform mesh comparing the strong and weak variational multiscale formulations (a) $L^2$ norm of velocity (b) $H^1_0$-semi norm of pressure.}
	\label{fig:Convergence}
\end{figure}

Another problem evaluated is the body force driven cavity taken from \cite{Donea}.  The problem geometry is the same as the Lid-driven cavity except that a velocity $v_x = v_y = 0.0$ is prescribed on the boundary and a constant body force is applied to the entire domain.  The prescribed constant body force is given as
\begin{align}
b_1 &= (12-24y)x^4 + (-24+48y)x^3 + (-48y+72y^2-48y^3+12)x^2 \nonumber \\
&+(-2+24y-72y^2+48y^3)x +1 - 4y +12y^2-8y^3 \nonumber \\
b_2 &= (8-48y+48y^2)x^3 + (-12+72y-72y^2)x^2 \nonumber \\ 
&+(4-24y+48y^2-48y^3+24y^4)x - 12y^2+24y^3-12y^4
\end{align}
The exact solution is
\begin{align}
v_x & = x^2(1-x)^2(2y-6y^2+4y^3) \nonumber \\
v_y & = -y^2(1-y)^2(2x-6x^2+4x^3) \nonumber \\
p & = x(1-x)
\end{align}
Numerical results are shown in Figure \ref{fig:BodyForce}, and they correspond well with other published results.  The convergence properties of the strong variational multiscale formulation are shown in Figure \ref{fig:Convergence}.  To measure the error in the velocity, the $L^2$ norm is used, whereas the $H^1$-semi norm is used to compute the error in the pressure.  Notice that the convergence rates are as expected for the Stokes problem using linear elements \cite{PSPG}.
%

%
%
%
%

\section{CONCLUSIONS}
We have introduced a new stabilized formulation for the Stokes problem that is appropriate for equal order interpolation for the velocity and pressure fields.  The new formulation produces a scalar stabilization parameter that is consistently derived purely by the method of weighted residuals. We have also shown that an equivalence between enriched finite element methods and stabilized methods for the Stokes problem does not exist for certain elements.  In particular, we have shown that enriching Q4 and B8 elements with standard bubble functions produces unstable results.  Clearly, this work highlights the need for more emphasis in the development of bubble function enriched methods and the exact nature of their relationship to stabilized formulations.

%
\section*{APPENDIX}
\subsection*{$\bullet$ Notation and definitions}
Consider an $n \ \times \ m$ matrix $\boldsymbol{A}$ and a $p \ \times \ q$ matrix $\boldsymbol{B}$
\begin{align*}
\boldsymbol{A} = \left[ \begin{array}{ccc} a_{1,1} & \hdots & a_{1,m}  \\
\vdots & \ddots & \vdots \\
a_{n,1} & \hdots & a_{n,m} \end{array} \right] \ ; \ \boldsymbol{B} = \left[ \begin{array}{ccc} b_{1,1} & \hdots & b_{1,q}  \\
\vdots & \ddots & \vdots \\
b_{p,1} & \hdots & b_{p,q} \end{array} \right]
\end{align*}
The \emph{Kronecker product} of these matrices is an $np \ \times \ mq$ matrix, and is defined as
\begin{align*}
\boldsymbol{A} \odot \boldsymbol{B} := \left[ \begin{array}{ccc} a_{1,1}\boldsymbol{B} & \hdots & a_{1,m}\boldsymbol{B}  \\
\vdots & \ddots & \vdots \\
a_{n,1}\boldsymbol{B} & \hdots & a_{n,m}\boldsymbol{B} \end{array} \right]
\end{align*}
The $\mathrm{vec}[\cdot]$ operator is defined as
\begin{align*}
\mathrm{vec}[\boldsymbol{A}] := \left[ \begin{array}{c} a_{1,1} \\
\vdots \\
a_{1,m} \\
\vdots \\
a_{n,1} \\
\vdots \\
a_{n,m} \end{array} \right]
\end{align*}
\subsection*{$\bullet$ The divergence of the jacobian matrix}
Consider the jacobian matrix  $\boldsymbol{J} := \frac{\partial \boldsymbol{x}}{\partial \boldsymbol{\xi}}$, a matrix of the coordinates of the nodes of an element, $\hat{\boldsymbol{x}}$, and the first and second derivatives of the shape functions $\boldsymbol{DN}$ and $\boldsymbol{D}^2\boldsymbol{N}$, such that $\boldsymbol{x} = \boldsymbol{N}\hat{\boldsymbol{x}}$, $[\boldsymbol{DN}]_{nm} = \frac{\partial\boldsymbol{N}_n}{\partial\boldsymbol{\xi}_m}$, and $[\boldsymbol{D}^2\boldsymbol{N}]_{nms} = \frac{\partial\boldsymbol{N}_n}{\partial\boldsymbol{\xi}_s\partial\boldsymbol{\xi}_m}$.  Starting with a simple identity, we can derive the divergence of the jacobian matrix as follows
\begin{align*}
  \boldsymbol{J} \boldsymbol{J}^{-1} & =  \boldsymbol{I}  \\
  J_{im} J^{-1}_{mk} & = \delta_{ik}  \\
  \frac{\partial}{\partial x_k} \left( J_{im} J^{-1}_{mk}  \right) & =  \frac{\partial}{\partial x_k} \left( \delta_{ik} \right) = 0 \\
  \frac{\partial J_{im}}{\partial x_k} J^{-1}_{mk} + J_{im} \frac{\partial J^{-1}_{mk}}{\partial x_k} & =  0  \\
  J^{-1}_{pi} \frac{\partial J_{im}}{\partial x_k} J^{-1}_{mk} + \delta_{pm} \frac{\partial J^{-1}_{mk}}{\partial x_k} & =  0  \\
\frac{\partial J^{-1}_{pk} }{\partial x_k} & =  -J^{-1}_{pi} \frac{\partial J_{im}}{\partial x_k} J^{-1}_{mk}  \\
\frac{\partial J^{-1}_{pk} }{\partial x_k} & =  -J^{-1}_{pi} \hat{x}_{ni} \frac{\partial }{\partial x_k}  \left(\frac{\partial N_n}{\partial \xi_m}\right)J^{-1}_{mk} \\
\frac{\partial J^{-1}_{pk} }{\partial x_k} & =  -J^{-1}_{pi} \hat{x}_{ni}\frac{\partial }{\partial \xi_s} \frac{\partial N_n}{\partial \xi_m} \frac{\partial \xi_s}{\partial x_k}J^{-1}_{mk}  \\
\frac{\partial J^{-1}_{pk} }{\partial x_k} = \frac{\partial}{\partial x_k} \frac{\partial \xi_p}{\partial x_k} & =  -J^{-1}_{pi}\hat{x}_{ni} D^2N_{nms} J^{-1}_{mk} J^{-1}_{sk}  \\
\nabla \cdot \boldsymbol{J}^{-1} & =  - \boldsymbol{J}^{-1} \hat{\boldsymbol{x}}^T \boldsymbol{D}^2\boldsymbol{N} \boldsymbol{J}^{-1} \boldsymbol{J}^{-T} 
\end{align*}
To further clarify, for a Q4 element, $\hat{\boldsymbol{x}}$, $\boldsymbol{DN}$, and $\boldsymbol{D}^2\boldsymbol{N}$ are defined as
\begin{align*}
  \hat{\boldsymbol{x}} &:= \left[ \begin{array}{cc} x_1 & y_1 \\
  x_2 & y_2 \\
  x_3 & y_3 \\
  x_4 & y_4 \\
 \end{array} \right] \\
 \boldsymbol{DN} &:= \left[ \begin{array}{cc} \frac{\partial N_1}{\partial \xi_1} & \frac{\partial N_1}{\partial \xi_2} \\
\vdots & \vdots \\
\frac{\partial N_4}{\partial \xi_1} & \frac{\partial N_4}{\partial \xi_2} \end{array} \right] \\
\boldsymbol{D}^2\boldsymbol{N} &:= \left[ \begin{array}{cccc} \frac{\partial^2 N_1}{\partial \xi_1 \partial \xi_1} & \frac{\partial^2 N_1}{\partial \xi_1 \partial \xi_2} & \hdots & \frac{\partial^2 N_1}{\partial \xi_2 \partial \xi_2} \\
\vdots & \vdots  & \ddots & \vdots  \\
\frac{\partial^2 N_4}{\partial \xi_1 \partial \xi_1} & \frac{\partial^2 N_4}{\partial \xi_1 \partial \xi_2} & \hdots & \frac{\partial^2 N_4}{\partial \xi_2 \partial \xi_2} \end{array} \right]
\end{align*}
\subsection*{$\bullet$ The Laplacian of a bubble function}
Noting that $\boldsymbol{\beta}$ is a constant vector and making use of th divergence of the jacobian matrix as shown above, the Laplacian of a bubble function can be computed as follows
\begin{align*}
\nabla^2 \boldsymbol{v}'  & =   \frac{\partial}{\partial x_i} \frac{\partial }{\partial x_i} b^e(\xi)\beta_j   \\
  & =   \frac{\partial}{\partial x_i} \left( \frac{\partial b^e(\xi)}{\partial \xi_k} \frac{\partial \xi_k}{\partial x_i} \right) \beta_j  \\
  & =   \left( \frac{\partial}{\partial x_i} \frac{\partial b^e(\xi)}{\partial \xi_k} \frac{\partial \xi_k}{\partial x_i} + \frac{\partial b^e(\xi)}{\partial \xi_k} \frac{\partial}{\partial x_i} \frac{\partial \xi_k}{\partial x_i} \right) \beta_j  \\
  & =   \left( \frac{\partial}{\partial \xi_p}\frac{\partial b^e(\xi)}{\partial \xi_k}\frac{\partial \xi_p}{\partial x_i}\frac{\partial \xi_k}{\partial x_i} + \frac{\partial b^e(\xi)}{\partial \xi_k} \frac{\partial}{\partial x_i} \frac{\partial \xi_k}{\partial x_i}\right)\beta_j   \\
  & =   \left( \frac{\partial}{\partial \xi_p}\frac{\partial b^e(\xi)}{\partial \xi_k}\frac{\partial \xi_p}{\partial x_i}\frac{\partial \xi_k}{\partial x_i} - \frac{\partial b^e(\xi)}{\partial \xi_k}\frac{\partial \xi_k}{\partial x_r}\hat{x}_{nr}\frac{\partial}{\partial \xi_s}\frac{\partial N_n}{\partial \xi_m}\frac{\partial \xi_m}{\partial x_i}\frac{\partial \xi_s}{\partial \xi_i} \right) \beta_j  \\
  & =   \left( H^{b^e}_{pk} J^{-1}_{pi} J^{-1}_{ki} - (\nabla_{\xi} b^e)_{k} J^{-1}_{kr} \hat{x}_{nr} D^2N_{nms} J^{-1}_{mi} J^{-1}_{si} \right) \beta_j  \\
  & =   \left( \boldsymbol{H}^{b^e} : \boldsymbol{J}^{-1} \boldsymbol{J}^{-T} - \nabla_{\xi} b^{eT} \boldsymbol{J}^{-1} \hat{\boldsymbol{x}}^T \boldsymbol{D}^2\boldsymbol{N}: \boldsymbol{J}^{-1} \boldsymbol{J}^{-T} \right)  \boldsymbol{\beta}  \\
  \nabla^2 b^e &:=  \left( \boldsymbol{H}^{b^e} : \boldsymbol{J}^{-1} \boldsymbol{J}^{-T} - \nabla_{\xi} b^{eT} \boldsymbol{J}^{-1} \hat{\boldsymbol{x}}^T \boldsymbol{D}^2\boldsymbol{N} : \boldsymbol{J}^{-1} \boldsymbol{J}^{-T} \right) 
\end{align*}
For trilinear B8 elements, $\boldsymbol{D}^2\boldsymbol{N}$ is defined as
\begin{align*}
\boldsymbol{D}^2\boldsymbol{N} &:= \left[ \begin{array}{ccccc} \frac{\partial^2 N_1}{\partial \xi_1 \partial \xi_1} & \frac{\partial^2 N_1}{\partial \xi_1 \partial \xi_2} & \frac{\partial^2 N_1}{\partial \xi_1 \partial \xi_3}& \hdots & \frac{\partial^2 N_1}{\partial \xi_3 \partial \xi_3} \\
\vdots & \vdots & \vdots  & \ddots & \vdots  \\
\frac{\partial^2 N_8}{\partial \xi_1 \partial \xi_1} & \frac{\partial^2 N_8}{\partial \xi_1 \partial \xi_2} & \frac{\partial^2 N_8}{\partial \xi_1 \partial \xi_3} & \hdots & \frac{\partial^2 N_8}{\partial \xi_3 \partial \xi_3}  \\
 \end{array} \right] 
\end{align*}
Note that $\boldsymbol{D}^2\boldsymbol{N}$ is a matrix representation of a third--order tensor.  The matrix representing the second derivatives of the bubble functions, $\boldsymbol{H}^{b^e}$, is given as 
\begin{align*}
\boldsymbol{H}^{b^e} &:= \left[ \begin{array}{ccccc} \frac{\partial^2 b^e}{\partial \xi_1 \partial \xi_1} & \frac{\partial^2 b^e}{\partial \xi_1 \partial \xi_2} & \frac{\partial^2 b^e}{\partial \xi_1 \partial \xi_3}& \hdots & \frac{\partial^2 b^e}{\partial \xi_3 \partial \xi_3} \end{array} \right] 
\end{align*}
%

\section*{ACKNOWLEDGMENTS}
The authors would like to thank Professor Arif Masud and Professor Albert Valocchi for their valuable insights provided in discussions on this topic. The research reported herein was partly 
supported (DZT) by the Computational Science and Engineering Fellowship at UIUC, and (KBN) by the Department of Energy through a SciDAC-2 project (Grant No. 
DOE DE-FC02-07ER64323). This support is gratefully acknowledged. The opinions 
expressed in this paper are those of the authors and do not necessarily reflect 
that of the sponsors. 

\bibliographystyle{unsrt}
\bibliography{../../DISSERTATION/References/references}
\newpage

%
\end{document}